\documentclass{amsart}
\usepackage{pstricks,pst-node,amssymb}
\usepackage{amssymb,psfig}
\usepackage{epsfig}
\usepackage{pstcol}             

\title{A Combinatorial Model for Crystals of Kac-Moody Algebras}

\author{Cristian Lenart and Alexander Postnikov}

\address{Department of Mathematics and Statistics, State University of New York at Albany, Albany, NY 12222}
\email{lenart@albany.edu}

\address{Department of Mathematics, M.I.T., Cambridge, MA 02139}
\email{apost@math.mit.edu}


\thanks{Cristian Lenart was supported by National Science Foundation 
grant DMS-0403029 and by SUNY Albany Faculty Research Award 1039703}
\thanks{Alexander Postnikov was supported by National Science Foundation 
grant DMS-0201494 and by Alfred P.\ Sloan Foundation research fellowship}


\date{February~4, 2005; updated June~4, 2005}

\numberwithin{equation}{section}

\theoremstyle{plain}
\newtheorem{theorem}{Theorem}[section]
\newtheorem{proposition}[theorem]{Proposition}
\newtheorem{lemma}[theorem]{Lemma}
\newtheorem{corollary}[theorem]{Corollary}
\newtheorem{conjecture}[theorem]{Conjecture}

\theoremstyle{definition}
\newtheorem{definition}[theorem]{Definition}
\newtheorem{example}[theorem]{Example}

\theoremstyle{remark}
\newtheorem{remark}[theorem]{Remark}
\newtheorem{remarks}[theorem]{Remarks}

\def\R{\mathbb{R}}
\def\Z{\mathbb{Z}}

\def\Q{\mathbb{Q}}

\def\Waff{W_{\mathrm{aff}}}

\def\hvee{h}
\def\h{\mathfrak{h}}
\def\hR{\mathfrak{h}^*_\mathbb{R}}
\newcommand{\casethree}[6]{\left\{ \begin{array}{ll} #1 &\mbox{if $#2$} \\#3 &\mbox{if $#4$} \\ #5 & \mbox{if $#6$}\,. \end{array} \right.}
\newcommand{\casethreeother}[5]{\left\{ \begin{array}{ll} #1 &\mbox{if $#2$} \\#3 &\mbox{if $#4$} \\ #5 & \mbox{otherwise}\,. \end{array} \right.}
\newcommand{\casetwo}[3]{\left\{ \begin{array}{ll} #1 &\mbox{if $#2$} \\ #3 &\mbox{otherwise}\,. \end{array} \right.}
\newcommand{\casetwoc}[3]{\left\{ \begin{array}{ll} #1 &\mbox{if $#2$} \\ #3 &\mbox{otherwise}\,, \end{array} \right.}
\newcommand{\casetwoex}[4]{\left\{ \begin{array}{ll} #1 &\mbox{if $#2$} \\ #3 &\mbox{if $#4$} \,. \end{array} \right.}

\newcommand{\casefour}[8]{\left\{ \begin{array}{ll} #1 &\mbox{if $#2$} \\#3 &\mbox{if $#4$} \\ #5 & \mbox{if $#6$} \\  #7 &\mbox{if $#8$}\,, \end{array} \right.}
\newcommand{\stacksum}[2]{\sum_{\begin{array}{c}\vspace{-6mm}\;\\ \vspace{-1mm}\scriptstyle{#1}\\ \scriptstyle{#2}\end{array}} }
\begin{document}
\bibliographystyle{plain}

\begin{abstract}
We present a simple combinatorial model for the characters of the irreducible
integrable highest weight modules for complex symmetrizable Kac-Moody algebras.
This model can be viewed as a discrete counterpart to the Littelmann path
model.  
We describe
crystal graphs and give a Littlewood-Richardson rule for decomposing tensor
products of irreducible representations. The new model is based on the notion
of a $\lambda$-chain, which is a chain of positive roots defined by certain
interlacing conditions.
\end{abstract}

\maketitle

\tableofcontents

\section{Introduction}
\label{sec:intro}


We have recently given a combinatorial model for the irreducible characters of
a complex semisimple Lie group $G$, and, more generally, for the Demazure
characters \cite{LP}. This model was defined in the context of the equivariant
$K$-theory of the generalized flag variety $G/B$.  Our character formulas were
derived from a Chevalley-type formula in $K_T(G/B)$.  Our model was based on
enumerating certain saturated chains in the Bruhat order on the
corresponding Weyl group $W$. This enumeration is determined by an 
alcove path, which is a sequence of adjacent alcoves for the affine Weyl
group $\Waff$ of the Langland's dual group $G^\vee$.  Alcove paths correspond
to decompositions of elements in the affine Weyl group.  Our Chevalley-type
formula was was expressed in terms of a certain  $R$-matrix, that is, in terms
of a collection of operators satisfying the  Yang-Baxter equation. This
setup allowed us to easily explain the independence of our formulas from the
choice of an alcove path. 

There are other models for Chevalley-type formulas in $K_T(G/B)$ and for the
characters of the irreducible representations of $G$ with highest weight
$\lambda$, such as the Littelmann path model.
Littelmann~\cite{Li1,Li2,Li3} showed that the characters can be described by
counting certain continuous paths in $\hR$.  These paths are constructed
recursively, by starting with an initial one, and by applying certain root
operators. By making specific choices for the initial path, one can obtain
special cases which have more explicit descriptions. For instance, a straight
line initial path leads to the Lakshmibai-Seshadri paths (LS paths).
These were introduced before Littelmann's work, in the context of standard
monomial theory~\cite{LS1}.  They have a nonrecursive description as weighted
chains in the Bruhat order on the quotient $W/W_\lambda$ of the corresponding
Weyl group $W$ modulo the stabilizer $W_\lambda$ of the weight $\lambda$;
therefore, we will use the term LS chains when referring to this
description.  LS paths were used by Pittie and Ram \cite{PR} to derive a
$K_T$-Chevalley formula. Recently, Gaussent and Littelmann~\cite{GaLi},
motivated by the study of Mirkovi\'c-Vilonen cycles, defined another
combinatorial model for the irreducible characters of a complex semisimple Lie
group.   This model is based on LS-galleries, which are certain sequences
of faces of alcoves for the corresponding affine Weyl group. According to
\cite{GaLi}, for each LS-gallery there is an associated Littelmann path and a
saturated chain in the Bruhat order on $W/W_\lambda$. 
Our alcove path model from~\cite{LP} and LS-galleries, which were developed
independently, can be related to each other in the case of regular weights.

The main goal of this paper is to extend the model introduced in \cite{LP} from
semisimple Lie algebras to complex symmetrizable Kac-Moody algebras. Instead of
alcove paths (that make sense only in finite types), the present paper is based
on $\lambda$-chains, which are chains of positive roots defined by
certain  interlacing conditions. These new objects extend the notion of a
reflection ordering \cite{Dyer}, and have many interesting combinatorial
properties, some of which are investigated in this paper. The extension from
the finite to the infinite case is nontrivial since the description in terms of
the affine Weyl groups and alcoves is no longer available in the latter case.
Compared to the geometric approach in \cite{LP} (involving the equivariant
$K$-theory of $G/B$), the generalization in this paper is carried out purely in
the context of representation theory, and, in this process, new features of our
model are added.  For instance, we define root operators in our model, and
study their properties.  This allows us to show that our model satisfies the
axioms of an admissible system of Stembridge \cite{St}.  Thus, we easily
derive character formulas, a Littlewood-Richardson rule for decomposing tensor
products of irreducible representations, as well as a branching rule. The
approach via admissible systems was already applied to LS chains in
\cite[Section 8]{St}.  Stembridge's approach has the advantage of making a part
of the proof independent of a particular model for Weyl characters, by using a
system of axioms for such models. 

The root operators in our model define certain colored graph structures which we call semiperfect crystals. In the case of a special $\lambda$-chain, we show that our model gives rise to LS paths (which are also defined in the Kac-Moody setup) via a certain process of passing to the limit. It then follows that the semiperfect crystals for  the special $\lambda$-chain coincide with Kashiwara's crystals (called here perfect crystals) of the corresponding irreducible integrable highest weight modules \cite{kascqa,kascbq}. Furthermore, the independence of our model from the choice of a $\lambda$-chain was proved in \cite{Le} in the finite case. Hence, in this case, we know that our semiperfect crystals are perfect crystals for any choice of a $\lambda$-chain. 

We briefly discuss the relationship between our model and the Littelmann path model (in the Kac-Moody setup). The former is a discrete counterpart of the latter without being simply a translation of it into a different language. Let us explain. Both models are based on certain choices (of a $\lambda$-chain and of an initial path, respectively). However, in its full generality, the path model is a recursive one, as explained above, whereas ours is a nonrecursive one. The definition of root operators in our model always involves reflecting a single segment in a $\lambda$-chain, whereas several segments of a path might be reflected in the path model (cf. Section \ref{sec:finite}). Furthermore, we have a generalization of the left-key construction of Lascoux-Sch\"utzenberger (relevant to Demazure characters, cf. Remark \ref{lskey}) in our model for any choice of a $\lambda$-chain, whereas, in the path model, a similar generalization exists only for LS paths. Although LS paths can be derived from our model, as explained above, the two are not equivalent either. Indeed, certain information is lost in the process of passing to the limit, namely the order in which certain hyperplanes are crossed, when a path is deformed into one which passes through the intersection of these hyperplanes. (In some cases, such as the one related to the $K_T$-Chevalley formula in \cite{PR}, the lost information needs to be recovered, but this is done by a nontrivial lift from $W/W_\lambda$ to $W$; the $K_T$-Chevalley formula in \cite{LP} based on our model obviates this lift.)

The main advantage of our model is related to its combinatorial nature. This model is computationally efficient  because it uses simple labels for the vertices of an irreducible crystal, namely certain finite sets of positive integers (called admissible subsets); these can be efficiently constructed. Although a $\lambda$-chain might be an infinite object, we do not need to store it, since we always work with a finite set of objects (cf. Remark \ref{finiteness}). More precisely, all we need is a procedure for comparing two elements in a $\lambda$-chain, namely two pairs, each consisting of a positive root and a nonnegative integer (see Proposition \ref{constr-lambdachain} for a simple procedure of this type). Let us also note that our model is equally
simple for both regular and nonregular highest weights $\lambda$. Indeed, instead of working with chains of cosets in $W/W_\lambda$ (as in the case of LS chains) or with possibly lower dimensional faces of alcoves (as in the case of LS-galleries), we always work with chains in $W$ and, in the finite case, with alcoves too. Compared to LS chains and LS-galleries, we also eliminate the need for making extra choices after the initial choice (of an initial path, of a sequence of alcoves/lower dimensional faces of them, or of a $\lambda$-chain). Unlike LS-galleries, our model extends to the Kac-Moody case.  Finally, it leads to an extensive generalization of the combinatorics of irreducible characters from Lie type $A$ (where the combinatorics is based on Young tableaux) to arbitrary type.  For example, an alcove path analog of Sch\"utzenberger's evacuation is given in \cite{Le}.  The study of the combinatorics of our model will be
continued in future publications. 

\medskip
Let us now present our combinatorial formula for characters.
Fix a complex symmetrizable Kac-Moody algebra $\mathfrak{g}$; 
see~\cite{kacidl}.
Let $\Phi^+$ be the associated set of positive real roots.
For a root $\alpha\in\Phi^+$, let $\alpha^\vee:=2\alpha/\langle \alpha,\alpha\rangle$
be the corresponding coroot.  Let $\lambda$ be a dominant integral weight, 
and let $V_\lambda$ be the associated irreducible representation of 
$\mathfrak{g}$ with the highest weight $\lambda$.

We define a  $\lambda$-chain (of roots) $\{\beta_i\}_{i\in I}$
as a map $I\to\Phi^+$, $i\mapsto \beta_i$, from a totally
ordered index set $I$ to positive roots, that satisfies the conditions below; given $i\in I$ and $\alpha\in\Phi^+$, we use the notation $N_i(\alpha):=\#\{j<i\mid \beta_j=\alpha\}$ and $N_i(-\alpha):=1-N_i(\alpha)$.
\begin{enumerate}
\item For a root $\alpha\in\Phi^+$, the number 
$\#\{i\in I\mid \beta_i=\alpha\}$
of occurrences of $\alpha$ equals $\langle \lambda,\alpha^\vee \rangle$.
\item Given any triple of roots $\alpha,\,\beta,\,\gamma$ such that $\alpha\ne\beta$ and $\gamma^\vee=\alpha^\vee+k\beta^\vee$ for some integer $k$, as well as $i\in I$ such that $\beta_i=\beta$, we have
\[N_i(\gamma)=N_i(\alpha)+k N_i(\beta)\,.\]
\end{enumerate}

It turns out that, in the finite case, condition~(2) is equivalent to a condition stating that any triple of roots $\alpha,\,\beta,\,\gamma$ with $\gamma^\vee=\alpha^\vee+\beta^\vee$ satisfies the following  interlacing property: 
there is exactly one element from the set $\{\alpha,\beta\}$ between
any two consecutive $\gamma$'s, as well as before the first $\gamma$; 
and there are no $\alpha$'s or $\beta$'s after the last $\gamma$.
Note that, according to condition~(1), the index set $I$ is always a
countable or a finite set.  However, it is not always order-isomorphic to 
a subset of $\Z$.  For example, it may contain infinite intervals.

\begin{lemma}\label{exist-lc}  For any dominant integral weight $\lambda$, 
there exists a $\lambda$-chain.
\end{lemma}

For a weight $\lambda$, there are usually many $\lambda$-chains.
We will give an explicit construction of a $\lambda$-chain (in Proposition \ref{constr-lambdachain}).

For $\alpha\in\Phi^+$, let $s_\alpha$ denote the associated reflection.
The reflections $s_\alpha$ generate the Weyl group $W$.
The covering relations $w \lessdot w s_\alpha$, for 
$w\in W$, $\ell(ws_\alpha)=\ell(w)+1$,
define the Bruhat order on $W$.
These covering relations are labelled by the corresponding positive roots $\alpha$.

Let us fix a $\lambda$-chain $\{\beta_i\}_{i\in I}$.  
Let us say that a finite subset $J=\{j_1<\dots<j_l\}$ of the index set 
$I$ is  admissible if the roots $\beta_{j_1},\dots,\beta_{j_l}$
are labels of an increasing saturated chain in the Bruhat order 
starting at the identity element, i.e., we have:
$$
1\lessdot s_{\beta_{j_1}}\lessdot s_{\beta_{j_1}} s_{\beta_{j_2}}\lessdot 
s_{\beta_{j_1}} s_{\beta_{j_2}} s_{\beta_{j_3}}
\lessdot 
\dots \lessdot s_{\beta_{j_1}}\dots s_{\beta_{j_l}}.
$$

For $\alpha\in\Phi^+$ and $k\in\Z$, let $s_{\alpha,k}$ be
the affine reflection given by $s_{\alpha,k}:\mu\mapsto s_\alpha(\mu) +
k\alpha$.  For a $\lambda$-chain, let $(k_i)_{i\in I}$ be the associated
sequence of nonnegative integers defined by 
$k_i := \#\{j<i\mid \beta_j = \beta_i\}$.
Let us define the  weight $\mu(J)$ of an admissible subset 
$J=\{j_1<\cdots < j_k\}$ by 
$$
\mu(J) = s_{\beta_{j_1},k_{j_1}}\cdots 
s_{\beta_{j_l},k_{j_l}}(\lambda).
$$

\begin{theorem}  For a dominant integral weight $\lambda$
and any $\lambda$-chain $\{\beta_i\}_{i\in I}$, the character
$\chi(\lambda)$ of the irreducible representation $V_\lambda$ is given by
$$
\chi(\lambda) = \sum_J e^{\mu(J)},
$$
where the sum is over the admissible subsets $J$ of the index set $I$.
\label{th:char_formula}
\end{theorem}

\begin{example}  Let $\mathfrak{g}$ be the Lie algebra of type $A_2$.
Let us fix a choice of simple roots $\alpha_1,\alpha_2$, 
and let $\lambda=\omega_1$ be the first fundamental weight.  
In this case, there is only one $\lambda$-chain
$(\beta_1,\beta_2)=(\alpha_1,\alpha_1+\alpha_2)$ (assuming that $I=\{1,2\}$).
The index set $I$ has 3 admissible subsets: $\emptyset, \{1\}, \{1,2\}$.
The subset $\{2\}$ is not admissible because the reflection 
$s_{\alpha_1+\alpha_2}$ does not cover the identity element.
In this case, $(k_1,k_2) = (0,0)$.
Theorem~\ref{th:char_formula} gives the following expression
for the character of $V_{\omega_1}$:
$$
\chi(\omega_1) = e^{\omega_1} + e^{s_{\alpha_1}(\omega_1)} + 
e^{s_{\alpha_1} s_{\alpha_1+\alpha_2}(\omega_1)}.
$$
\end{example}

We will define a colored graph structure on the collection of
admissible subsets.  From this, we will deduce a rule for decomposing  
tensor products of irreducible representations and a branching rule.

\medskip
The general outline of the paper follows.  In Section~\ref{sec:notation},
we review basic notions related to roots systems for complex symmetrizable
Kac-Moody algebras and fix our notation.  In Section~\ref{sec:admissible}, 
we discuss crystals and give Stembridge's axioms of admissible systems.
In Sections~\ref{sec:lambdachains}-\ref{sec:folded}, we define our combinatorial
model.  In Section~\ref{sec:lambdachains}, we discuss $\lambda$-chains.
In Section~\ref{sec:chains-foldings}, we define and study 
folding operators; we use them to construct more general chains of 
roots from a $\lambda$-chain, which we call admissible foldings.  
In Section~\ref{sec:folded}, we study combinatorial properties 
of admissible subsets and admissible foldings. 
In Section~\ref{sec:rootoperators}, we define root
operators on admissible subsets/foldings.  In
Section~\ref{sec:mainthm}, we prove that our combinatorial model satisfies 
Stembridge's axioms.  This enables us to derive a character formula
for complex symmetrizable Kac-Moody algebras, a Littlewood-Richardson rule for
decomposing tensor products of irreducible representations, as well as a
branching rule. In Section \ref{sec:ls}, we explain the connection between our
model and LS chains. In Section \ref{sec:finite}, we discuss the way in which the present
model specializes to the one in our previous paper~\cite{LP} 
in the finite case.


\medskip
{\bf Acknowledgments.}  We are grateful to 
John Stembridge for the explanation of his work, and to
V.~Lakshmibai and Peter Magyar for
helpful conversations.

\section{Preliminaries}
\label{sec:notation}

In this section, we briefly recall the general setup for complex symmetrizable Kac-Moody algebras and their representations. We closely follow \cite[Section 1]{St}, and refer to \cite{kacidl,kumkmg} for more details. 

Let $V$ be a finite-dimensional real vector space with a nondegenerate symmetric bilinear form $\langle\,\cdot\,,\,\cdot\,\rangle$, and let $\Phi\subset V$ be a crystallographic root system of rank $r$ with simple roots $\{\alpha_1,\ldots,\alpha_r\}$. By this, we mean that $\Phi$ is the set of real roots of some complex symmetrizable Kac-Moody algebra. The finite root systems of this type are the root systems of semisimple Lie algebras. Note that, in the infinite case, it is possible for the simple roots to span a proper subspace of $V$; indeed, it can happen that the bilinear form is degenerate on the span of the simple roots. 

Given a root $\alpha$, the corresponding {\it coroot\/} is $\alpha^\vee := 2\alpha/\langle \alpha,\alpha \rangle$.  The collection of coroots
$\Phi^\vee:=\{\alpha^\vee \mid  \alpha\in\Phi\}$ forms the 
{\it dual root system.} For each root $\alpha$, there is a reflection $s_\alpha\::\: V\to V$ defined by
\[s_\alpha:  \lambda \mapsto \lambda - \langle \lambda,\alpha^\vee \rangle\,\alpha.\]
More generally, for any integer $k$, one can consider the affine hyperplane 
\[H_{\alpha,k}:=\{\lambda\in V\mid  \langle \lambda,\alpha^\vee \rangle=k\}\,,\]
and let $s_{\alpha,k}$ denote the corresponding reflection, that is, 
\begin{equation}\label{affref}s_{\alpha,k}\,:\,\lambda\mapsto s_\alpha(\lambda)+k\alpha\,.\end{equation}

The {\it Weyl group\/} $W$ is the subgroup of $GL(V)$ 
 generated by the reflections $s_{\alpha}$ for $\alpha\in\Phi$.
In fact, the Weyl group $W$ is a Coxeter group, which is generated by the
{\it simple reflections\/} $s_1,\dots,s_r$ corresponding 
to the simple roots $s_p := s_{\alpha_p}$, subject to the 
{\it Coxeter relations}:
$$(
s_p)^2=1
\quad\textrm{and}\quad 
(s_p s_q)^{m_{pq}}=1\,;
$$
here the relations of the second type correspond to the distinct $p,\,q$ in $\{1,\ldots,r\}$ for which the dihedral subgroup generated by $s_p$ and $s_q$ is finite, in which case $m_{pq}$ is half the order of this subgroup. The Weyl group is finite if and only if $\Phi$ is finite.

An expression of a Weyl group element $w$ as a product 
of generators $w=s_{p_1}\cdots s_{p_l}$
which has minimal length is called a {\it reduced decomposition\/}
for $w$; its length $\ell(w)=l$ is called the {\it length\/} of $w$.
For $u,w\in W$, we say that $u$ {\it covers\/} $w$, and write $u\gtrdot w$,
if $w=u s_{\beta}$, for some $\beta\in\Phi^+$, and $\ell(u)=\ell(w)+1$.
The transitive closure ``$>$'' of the relation ``$\gtrdot$'' is called
the {\it Bruhat order\/} on $W$.

Let us note that $\Phi$ can be characterized by the following axioms:
\begin{enumerate}
\item[(R1)] $\{\alpha_1,\ldots,\alpha_r\}$ is a linearly independent set.
\item[(R2)] $\langle\alpha_p,\alpha_p\rangle>0$ for all $p=1,\ldots,r$.
\item[(R3)] $\langle\alpha_p,\alpha_q^\vee\rangle\in{\mathbb Z}_{\le 0}$ for all distinct simple roots $\alpha_p$ and $\alpha_q$. 
\item[(R4)] $\Phi=\bigcup_{p=1}^r W\alpha_p$. 
\end{enumerate}

Let $\Phi^+\subset \Phi$ be the set 
of positive roots, that is, the set of roots in the nonnegative linear span of the simple roots.  
Then $\Phi$ is the disjoint union of $\Phi^+$ and $\Phi^- := -\Phi^+$. 
We write $\alpha>0$ (respectively, $\alpha<0$) for $\alpha\in\Phi^+$ (respectively, 
$\alpha\in\Phi^-$), and we define ${\rm sgn}(\alpha)$ to be $1$ (respectively, $-1$). 
We also use the notation $|\alpha|:={\rm sgn}(\alpha)\alpha$. 

Let us choose a set of {\it fundamental weights\/}
$\omega_1,\dots,\omega_r$, that is, $r$ vectors in V such that $\langle\omega_p,\alpha_q^\vee\rangle=\delta_{pq}$. The lattice of (integral) {\em weights} $\Lambda$ is the lattice generated by the fundamental weights. 
The set $\Lambda^+$ of {\it dominant weights\/} is given by
$$
\Lambda^+:=\{\lambda\in\Lambda \mid  \langle \lambda,\alpha^\vee \rangle\geq 0
\textrm{ for any } \alpha\in\Phi^+\}.
$$
If we replace the weak inequalities above with strict ones, we obtain the strongly dominant weights. It is known that every $W$-orbit in $V$ has at most one dominant member. 

The (integral) Tits cone $\Lambda_c$ is defined to be the union of all $W$-orbits of dominant integral weights, or, equivalently,
\[\Lambda_c:=\{\lambda\in\Lambda\mid \langle\lambda,\alpha^\vee\rangle<0 \textrm{ for finitely many } \alpha\in\Phi^+\}\,.\]
We have $\Lambda=\Lambda_c$ in the finite case, but not otherwise.

We now define a ring $R$ that contains the characters of all irreducible representations of the corresponding Kac-Moody algebra. In the finite case, one may simply take $R$ to be the group ring of $\Lambda$, but in general more care is required.

First, we choose a height function $\mathrm{ht}\::\:V\to {\mathbb R}$, that is, a linear map assigning the value 1 to all simple roots. Second, for each $\lambda\in\Lambda$, let $e^\lambda$ denote a formal exponential subject to the rules $e^\mu\cdot e^\nu=e^{\mu+\nu}$ for all $\mu,\nu\in\Lambda$. We now define the ring $R$ to consist of all formal sums $\sum_{\lambda\in\Lambda}c_\lambda e^\lambda$ with $c_\lambda\in{\mathbb Z}$ satisfying the condition that there are only finitely many weights $\lambda$ with $\mathrm{ht}(\lambda)>h$ and $c_\lambda\ne 0$, for all $h\in{\mathbb R}$. The ring $R$ contains the formal power series ring $R_0={\mathbb Z}[[e^{-\alpha_1},\ldots,e^{-\alpha_r}]]$. In particular, if $f\in R_0$ has constant term 1, then $e^\lambda f$ has a multiplicative inverse in $R$.

For each $\lambda\in\Lambda^+$ with a finite $W$-stabilizer, we define 
\[\varDelta(\lambda):=\sum_{w\in W}\mathrm{sgn}(w)e^{w(\lambda)}\,,\]
where $\mathrm{sgn}(w)=(-1)^{\ell(w)}$. It is not hard to check that $\varDelta(\lambda)$ is a well-defined member of $R$. Moreover, if $\lambda$ is dominant, then $\varDelta(\lambda)\ne 0$ if and only if $\lambda$ is strongly dominant. In that case, $e^{-\lambda} \varDelta(\lambda)\in R_0$ has constant term 1, and $\varDelta(\lambda)$ is invertible in $R$.

Since the scalar product is nondegenerate, we may select $\rho\in\Lambda^+$ so that $\langle\rho,\alpha_p^\vee\rangle=1$ for all $p=1,\ldots,r$. This given, for each $\lambda\in\Lambda^+$ we define
\[\chi(\lambda):=\frac{\varDelta(\lambda+\rho)}{\varDelta(\rho)}=\frac{\sum_{w\in W}\mathrm{sgn}(w)e^{w(\lambda+\rho)-\rho}}{\sum_{w\in W}\mathrm{sgn}(w)e^{w(\rho)-\rho}}\in R\,.\]
It is easy to show that $w(\rho)-\rho$, and hence $\chi(\lambda)$, do not depend on the choice of $\rho$. By the Kac-Weyl character formula \cite{kacidl}, $\chi(\lambda)$ is the character of the irreducible integrable module of highest weight $\lambda$ for the corresponding Kac-Moody algebra. 

\section{Crystals}
\label{sec:admissible}

This section closely follows \cite[Section 2]{St}. We refer to this paper for more details.

\begin{definition}(cf.\ \cite{St}). A {\em crystal} is a 4-tuple $(X,\mu,\delta,\{F_1,\ldots,F_r\})$ satisfying Axioms (A1)-(A3) below, where
\begin{itemize}
\item $X$ is a set whose elements are called objects;
\item $\mu$ and $\delta$ are maps $X\rightarrow\Lambda$;
\item $F_p$ are bijections between two subsets of $X$.
\end{itemize}
A crystal is called an {\em admissible system} if it satisfies Axioms (A0) and (A4). A crystal is called {\em semiperfect} if it satisfies Axioms (A4) and (A5).
\end{definition}

\vspace{2mm}
{\bf (A0)} For all real numbers $h$, there are only finitely many objects $x$ such that $\mathrm{ht}(\mu(x))>h$.
\vspace{2mm}

Axiom (A0) implies that the generating series $G_X:=\sum_{x\in X}e^{\mu(x)}$ is a well-defined member of $R$. 

For each $x\in X$, we call $\mu(x)$, $\delta(x)$, and $\varepsilon(x):=\mu(x)-\delta(x)$ the {\em weight}, {\em depth}, and {\em rise} of $x$.

\vspace{2mm}
{\bf (A1)} $\delta(x)\in-\Lambda^+$, $\varepsilon(x)\in\Lambda^+$.
\vspace{2mm}

We define the depth and rise in the direction $\alpha_p$ by $\delta(x,p):=\langle \delta(x),\alpha_p^\vee\rangle$ and $\varepsilon(x,p):=\langle \varepsilon(x),\alpha_p^\vee\rangle$. In fact, we will develop the whole theory in terms of $\delta(x,p)$ and $\varepsilon(x,p)$ rather than $\delta(x)$ and $\varepsilon(x)$. 

\vspace{2mm}
{\bf (A2)} $F_p$ is a bijection from $\{x\in X\mid \varepsilon(x,p)>0\}$ to $\{x\in X\mid \delta(x,p)<0\}$.
\vspace{2mm}

We let $E_p:=F_p^{-1}$ denote the inverse map.

\vspace{2mm}
{\bf (A3)} $\mu(F_p(x))=\mu(x)-\alpha_p$, $\delta(F_p(x),p)=\delta(x,p)-1$.
\vspace{2mm}

Hence, we also have $\varepsilon(F_p(x),p)=\varepsilon(x,p)-1$.  The maps $E_p$ and $F_p$ act as raising and lowering operators that provide a partition of the objects into $\alpha_p$-strings that are closed under the action of $E_p$ and $F_p$. For example, the $\alpha_p$-{\em string through} $x$ is (by definition)
\[F_p^\varepsilon(x),\,\ldots,\,F_p(x),\,x,\,E_p(x),\,\ldots,\,E_p^{-\delta}(x)\,,\]
where $\delta=\delta(x,p)$ and $\varepsilon=\varepsilon(x,p)$. 

We define partial orders on $X$, one for each $p$, by
\begin{equation}\label{partialord}x\preceq_py\;\;\;\mathrm{if}\;\;\;x=F_p^k(y)\;\:\textrm{for some }\,k\ge 0\,.\end{equation}
Any assignment of elements $t(x,p)$ of a totally ordered set to pairs $(x,p)$ with $\delta(x,p)<0$ is called a {\em timing pattern} for $X$. A timing pattern is called {\em coherent} if the following two conditions are satisfied for all pairs $(x,p)$ such that $\delta(x,p)<0$ and $\varepsilon(x,p)>0$:
\begin{itemize}
\item $t(x,p)\ge t(F_p(x),p)$;
\item for all $q\ne p$, all integers $\delta<0$, and all $t\ge t(x,p)$, there is an object $y\succeq_q x$ such that $\delta(y,q)=\delta$ and $t(y,q)=t$ if and only if there is an object $y'\succeq_q F_p(x)$ such that $\delta(y',q)=\delta$ and $t(y',q)=t$. 
\end{itemize}
The following axiom ensures the existence of a certain sign-reversing involution used to cancel the negative terms in the Kac-Weyl character formula. 

\vspace{2mm}
{\bf (A4)} There exists a coherent timing pattern for $X$. 
\vspace{2mm}

Note that, compared to \cite{St}, here we let the timing pattern take values in any totally ordered set, and we reverse the total order previously considered. However, these minor changes, dictated by our needs, are easily taken care of by minor changes in the corresponding proofs in \cite{St}.  

We say that $X$ has a maximum object if it does so with respect to the partial order generated by all partial orders $\preceq_p$, for $p=1,\ldots,r$. 

\vspace{2mm}
{\bf (A5)} $X$ has a maximum object. 
\vspace{2mm}

\begin{proposition} Semiperfect crystals are admissible systems. \end{proposition}

\begin{proof} It suffices to show that Axiom (A0) is a consequence of the other axioms. Let $x_0$ be the maximum object of $X$ and let $\lambda:=\mu(x_0)$. Any
other object of $X$ is obtained from $x_0$ by applying lowering operators $F_p$. Each of these operators decreases weight by $\alpha_p$, and thus they decrease height by 1.  It follows that we may have at most $r^k$ elements of height 
${\rm ht}(\lambda)-k$, where $r$ is the number of simple roots.
\end{proof}

\begin{theorem}\cite{St}\label{thmst} If $X$ is an admissible system, then 
\[\chi(\nu)\cdot G_X=\sum_{x\in X\::\:\nu+\delta(x)\in\Lambda^+}\chi(\nu+\mu(x))\,.\]
In particular, if $X$ is a semiperfect crystal with maximal object $x_*$, then $G_X=\chi(\mu(x_*))$.
\end{theorem}

Given $P\subseteq\{1,\ldots,r\}$, let $\Phi_P$ denote the root subsystem of $\Phi$ with simple roots $\{\alpha_p\mid p\in P\}$. Following \cite{St}, we let $W_P\subseteq W$, $\Lambda_P\supseteq\Lambda$, and $R_P$ denote the corresponding Weyl group, weight lattice, and character ring. Provided that we use the height function inherited from $\Phi$ (in which case $R_P\supseteq R$), it is easy to see that any admissible system $X$ can also be viewed as an admissible system relative to $\Phi_P$ using only the operators $E_p$ and $F_p$ for $p\in P$. Given $\lambda\in\Lambda_P^+$, we let $\chi(\lambda;P)\in R_P$ denote the Weyl character (relative to $\Phi_P$) corresponding to $\lambda$. The following branching rule is given in \cite{St}.

\begin{corollary}\cite{St}\label{branching} If $X$ is an admissible system and $P\subseteq\{1,\ldots,r\}$, we have 
\[G_X=\sum_{x\::\:\delta(x)\in\Lambda_P^+}\chi(\mu(x);P)\,.\]
\end{corollary}

Finally, note that one can define on $X$ the structure of a directed colored graph by constructing arrows $x\rightarrow y$ colored $p$ for each $F_p(x)=y$. Also note that the most important source of crystals are integrable representations of quantum groups $U_q({\mathfrak g})$ corresponding to Kac-Moody algebras ${\mathfrak g}$. It is well-known that such representations can be encoded by combinatorial data $(X,\mu,\delta,\{F_1,\ldots,F_r\})$ satisfying Axioms (A1)-(A3), i.e., by crystals (cf. \cite{kascqa,kascbq}, see also e.g. \cite{josqgp}). These crystals are usually called {\em perfect crystals}, while the corresponding directed colored graphs are known as {\em crystal graphs}.

\section{$\lambda$-Chains of Roots}\label{sec:lambdachains}

Fix a dominant weight $\lambda$. Throughout this paper, we will use the term ``sequence'' for any map $i\mapsto a_i$ from a totally ordered set $I$ to some other set. We will use the notation $\{a_i\}_{i\in I}$, and, if $I$ is finite or countable, the usual notation $(a_{i_1},\,a_{i_2},\ldots)$, where $I=\{i_1<i_2<\ldots\}$. Given an element $b$ and an index $i$, we also define the following counting functions:
\[N(b):=\#\{k\mid  a_k=b\},\; N_i(b):=\#\{k<i\mid  a_k=b\},\; N_{ij}(b):=\#\{i<k<j\mid  a_k=b\},\]
assuming that the corresponding cardinalities are finite. If the elements of a sequence are positive roots and $\alpha$ is such a root, then we set $N_i(-\alpha):=1-N_i(\alpha)$. 

\begin{definition}\label{lambdachain} A $\lambda$-{\em chain (of roots)} is a sequence 
of positive roots $\{\beta_i\}_{i\in I}$ indexed by the elements of a totally ordered 
set $I$, which satisfies the conditions below. 
\begin{enumerate}
\item The number of occurrences of any positive root $\alpha$, i.e., $N(\alpha)$, is 
$\langle \lambda,\alpha^\vee \rangle$. 
\item Given any triple of roots $\alpha,\,\beta,\,\gamma$ such that $\alpha\ne\beta$ and $\gamma^\vee=\alpha^\vee+k\beta^\vee$ for some integer $k$, as well as $i\in I$ such that $\beta_i=\beta$, we have
\[N_i(\gamma)=N_i(\alpha)+k N_i(\beta)\,.\]
\end{enumerate}
\end{definition}


Note that finding a $\lambda$-chain amounts to defining a total order on the set 
\begin{equation}\label{indset}
I:=\{(\alpha,k)\mid  \alpha\in\Phi^+,\,0\le k<\langle \lambda,\alpha^\vee \rangle\}
\end{equation}
 such that condition (2) above holds, where $\beta_i=\alpha$ for any
$i=(\alpha,k)$ in $I$.  One particular example of such an order can be constructed
as follows. Fix a total
order on the set of simple roots $\alpha_1<\alpha_2<\ldots<\alpha_r$. For each
$i=(\alpha,k)$ in $I$, let
$\alpha^\vee=c_1\alpha_1^\vee+\ldots+c_r\alpha_r^\vee$, and define the vector 
\begin{equation}
v_i:=\frac{1}{\langle \lambda,\alpha^\vee \rangle}(k,c_1,\ldots,c_r)
\end{equation}
in $\Q^{r+1}$. 
The map $i\mapsto v_i$ is injective. Indeed, assume that $v_i=v_{i'}$ for $i=(\alpha,k)$ and $i'=(\alpha',k')$. If $\alpha\ne\alpha'$,  
the root system $\Phi^\vee$ would contain two proportional positive coroots 
$\alpha^\vee\ne(\alpha')^\vee$, which is not possible.
Also, the fact that $\alpha=\alpha'$ implies that $k=k'$. Hence, we can define a total order on $I$ by $i<j$ iff $v_i<v_j$ in the lexicographic order on $\Q^{r+1}$. We are now ready to prove the existence of $\lambda$-chains (cf. Lemma~\ref{exist-lc}).

\begin{proposition}\label{constr-lambdachain} Given the total order on $I$ defined above, the sequence $\{\beta_i\}_{i\in I}$ defined by $\beta_i=\alpha$ for $i=(\alpha,k)$ is a $\lambda$-chain.
\end{proposition}

\begin{proof} Throughout this proof, we work with coroots only and modify the above notions accordingly. Let us fix a triple of positive coroots $\alpha,\,\beta,\,\gamma$ such that $\gamma=\alpha+k\beta$ for some fixed integer $k>0$. Let $\beta':=k\beta$ and 
\[\widehat{I}:=\{(\delta,l)\mid  \delta\in(\Phi^\vee)^+,\,0\le l<\langle \lambda,\delta \rangle\}\cup\{(\beta',m)\mid 0\le m<\langle\lambda,\beta'\rangle\}\,.\]
We define a map $i\mapsto v_i$ for $\widehat{I}$ to ${\mathbb Q}^{r+1}$ as above. This map is not injective, because we have $v_{(\beta,l)}=v_{(\beta',kl)}$ for $0\le l<\langle\lambda,\beta\rangle$. However, we can still define a total order $\widehat{I}$ if, in addition to setting $i<j$ whenever $v_i<v_j$ in the lexicographic order, we impose $(\beta,l)<(\beta',kl)$ for $0\le l<\langle\lambda,\beta\rangle$. Hence, we can define the sequence $\{\delta_i\}_{i\in \widehat{I}}$ by $\delta_i=\delta$ for $i=(\delta,l)$. 

The remainder of the proof consists in showing that the finite sequence $\{\delta_j\}_{j\in \widehat{J}}$, where $\widehat{J}:=\{j\in \widehat{I}\mid  \delta_j\in\{\alpha,\beta',\gamma\}\}$ has the induced total order, is a concatenation of pairs $(\alpha,\gamma)$ and $(\beta',\gamma)$ (in any order). Hence, given $i\in\widehat{I}$ such that $\delta_i=\beta$, we have $N_i(\gamma)=N_i(\alpha)+N_i(\beta')$. On the other hand, it is easy to see that $N_i(\beta')=k N_i(\beta)$. Thus, condition (2) in the definition of a $\lambda$-chain is verified in the case when $\alpha$ and $\gamma$ are both positive or both negative roots. The remaining case is checked in a similar way. 

For each $p$ in $\{1,\ldots,r\}$, let us denote by $c_p^\alpha$, $c_p^{\beta'}$, $c_p^\gamma$ the coefficients of $\alpha_p^\vee$ in $\alpha$, $\beta'$, $\gamma$, respectively; clearly, $c_p^\gamma=c_p^\alpha+c_p^{\beta'}$. Also, let $a:=\langle\lambda,\alpha\rangle$, $b:=\langle\lambda,\beta'\rangle$, and $c:=\langle\lambda,\gamma\rangle$, where $c=a+b$. The interlacing condition to be verified is straightforward if $a=0$ or $b=0$. Assume $a\ne 0\ne b$. In this case, the mentioned condition is checked based on the following three claims about the finite sequence $\{\delta_j\}_{j\in \widehat{J}}$:
\begin{enumerate}
\item if an entry $\alpha$ is followed by $\beta'$, or vice versa, there is an entry $\gamma$ in-between;
\item between two entries $\alpha$, or two entries $\beta'$, there is an entry $\gamma$;
\item the sequence cannot start with $\gamma$, but must end with $\gamma$.
\end{enumerate}

For the first claim, let $(\alpha,l)<(\beta',m)$, which means $l/a\le m/b$. If the inequality is strict, then the rational number $(l+m)/c=(l+m)/(a+b)$ is strictly in-between the previous two ones; therefore, we have $(\alpha,l)<(\gamma,l+m)<(\beta',m)$, which proves the claim. Otherwise, all three numbers are equal. We can now repeat our reasoning above with $l$ and $m$ replaced by $c_p^\alpha$ and $c_p^{\beta'}$ for $p=1$. If equality still holds, we let $p=2$ etc. At some point, we must have strict inequalities; indeed, otherwise $\alpha$ and $\beta'$ would be proportional, which is impossible. 

For the second claim, we can assume that we have two consecutive entries $\alpha$ with corresponding indices $(\alpha,l-1)$ and $(\alpha,l)$. Furthermore, we can assume that we have an entry $\beta'$ with corresponding index $(\beta',m)$ such that
\[\frac{m}{b}\le\frac{l-1}{a}<\frac{l}{a}\le\frac{m+1}{b}\le 1\,.\]
It is straightforward to check that we have
\[\frac{l-1}{a}<\frac{l+m}{a+b}<\frac{l}{a}\,;\]
therefore, we have $(\alpha,l-1)<(\gamma,l+m)<(\alpha,l)$, which proves the claim.

For the last claim, note that the first three entries in the sequence $\{\delta_j\}_{j\in \widehat{J}}$ are $\alpha,\,\beta',\,\gamma$ in some order, and that the corresponding vectors in $\Q^{r+1}$ have their first components equal to $0$. The fact that the sequence cannot start with $\gamma$ follows by an argument similar to the proof of the first claim above. On the other hand, the sequence must end with $\gamma$ because the largest value of the first component in the vectors involved is $(c-1)/c$, and this value only appears in $v_i$ for $i=(\gamma,c-1)$.
\end{proof}

The following immediate consequence of condition (2) in the definition of a $\lambda$-chain will be used several times below. 

\begin{lemma}\label{countroots}
For any two positive roots $\alpha\ne\beta$ and $i\in I$ such that $\beta_i=\beta$, we have 
\[ N_i(s_\beta(\alpha))=N_i(\alpha)-\langle \beta,\alpha^\vee\rangle N_i(\beta)\,.\]
\end{lemma}

We conclude this section by showing that, in the case of a finite root system, we can simplify the definition of a $\lambda$-chain.

\begin{proposition}\label{equivalence} 
In the case of a finite root system, condition {\rm (2)} in Definition~{\rm \ref{lambdachain}} is equivalent to the following interlacing condition.
\begin{enumerate}
\item[(2$'$)] For each triple of positive roots $(\alpha,\,\beta,\,\gamma)$ with $\gamma^\vee=\alpha^\vee+\beta^\vee$, the finite sequence $\{\beta_j\}_{j\in J}$, where $J:=\{j\in I\mid  \beta_j\in\{\alpha,\beta,\gamma\}\}$ has the induced total order, is a concatenation of pairs $(\alpha,\gamma)$ and $(\beta,\gamma)$ (in any order). 
\end{enumerate}
\end{proposition}

\begin{proof} It is easy to see that the interlacing condition (2$'$) is equivalent to the following one: given any $i\in I$ and any triple of positive roots $(\alpha,\,\beta,\,\gamma)$ with $\gamma^\vee=\alpha^\vee+\beta^\vee$, we have
\begin{equation}\label{basiccount}
N_i(\gamma)=\casetwoex{N_i(\alpha)+N_i(\beta)}{\beta_i=\alpha\,\textrm{ or }\,\beta_i=\beta}{N_i(\alpha)+N_i(\beta)-1}{\beta_i=\gamma}
\end{equation}
Clearly, condition (2) in Definition~\ref{lambdachain} implies (\ref{basiccount}); indeed, one just has to set $k=1$ in the former. 

Assume that $\alpha$, $\beta$, and $\gamma$ are all positive roots and $k\ge 0$. Given that (\ref{basiccount}) holds, we check condition (2) in this case using induction on $k$, which starts at $k=0$. Given $k>0$, we know that $\alpha^\vee+l\beta^\vee$, for $l=1,\ldots,k-1$, are coroots, as elements of the $\beta^\vee$-string through $\alpha^\vee$. The induction step consists in writing $\gamma^\vee=(\alpha^\vee+(k-1)\beta^\vee)+\beta^\vee$, and in combining (\ref{basiccount}) with the induction hypothesis. 

Now assume that $\beta$ and $\gamma$ are positive roots, but $\alpha$ is a negative root. We check condition (2) in this case in a similar way, using induction which starts at $k=1$. Given an arbitrary $k>1$, we know that $-\alpha^\vee+l\beta^\vee$, for $l=1,\ldots,k-1$, are coroots, as elements of the $\beta^\vee$-string through $-\alpha^\vee$. We now write 
\begin{align}\gamma^\vee&=(-\alpha^\vee+(k-1)\beta^\vee)+\beta^\vee\;\;\mbox{if}\;\;-\alpha^\vee+(k-1)\beta^\vee>0\,,\;\;\mbox{and}\label{eq1}\\
\beta^\vee&=\gamma^\vee+(\alpha^\vee-(k-1)\beta^\vee)\;\;\mbox{otherwise}\,.\label{eq2}
\end{align}
The induction step is completed by combining (\ref{basiccount}) with the induction hypothesis (in (\ref{eq1})) and condition (2) in the case already verified (in (\ref{eq2})).
\end{proof}


For the rest of our construction (Sections~\ref{sec:chains-foldings}--\ref{sec:mainthm}),
let us fix a dominant integral weight $\lambda$ and fix an arbitrary 
$\lambda$-chain $\{\beta_i\}_{i\in I}$. 
We will use the notation $r_i$ for the reflection $s_{\beta_i}$.

\section{Folding Chains of Roots}\label{sec:chains-foldings}

We start by associating to our fixed $\lambda$-chain the closely related object
\[\Gamma(\emptyset):=(\{(\beta_i,1)\}_{i\in I},\rho) \,,\]
where  $\rho$ is a fixed dominant weight satisfying $\langle\rho,\alpha_p^\vee\rangle=1$ for all $p=1,\ldots,r$. Here, as well as throughout this article, we let $\infty$ be greater than all elements in $I$. We use operators called {\em folding operators} to construct from $\Gamma(\emptyset)$ new objects of the form 
\begin{equation}\label{defgamma}\Gamma=(\{(\gamma_i,\varepsilon_i)\}_{i\in I},\gamma_\infty) \,;\end{equation}
here $\gamma_i$ is a root, $\varepsilon_i=\pm 1$, any given root appears only finitely many times in $\Gamma$, and $\gamma_\infty$ is in the $W$-orbit of $\rho$. 
More precisely, given $\Gamma$ as above and $i$ in $I$, we let $t_i:=s_{\gamma_i}$ and we define 
\[\phi_i(\Gamma):=(\{(\delta_j,\zeta_j)\}_{j\in I},{t}_i(\gamma_\infty)) \,,\]
where
\[(\delta_j,\zeta_j):=\casethree{(\gamma_j,\varepsilon_j)}{j<i}{(\gamma_j,-\varepsilon_j))}{j=i}{(t_i(\gamma_j),\varepsilon_j)}{j>i}\]

Let us now consider the set of all $\Gamma$ that are obtained from $\Gamma(\emptyset)$ by applying folding operators; we call these objects the {\em foldings} of $\Gamma(\emptyset)$. Clearly, $\phi_i$ is an involution on the set of foldings of $\Gamma(\emptyset)$. In order to describe this set, let us note that the folding operators commute. Indeed, if $\Gamma$ is as in (\ref{defgamma}), $i<j$, $\alpha:=\gamma_i$, and $\beta:=\gamma_j$, then we have 
\[t_it_j=s_\alpha s_\beta= (s_\alpha s_\beta s_\alpha)s_\alpha=s_{s_\alpha(\beta)}s_\alpha\,;\]
 so $\phi_i\phi_j(\Gamma)=\phi_j\phi_i(\Gamma)$. This means that every folding $\Gamma$ of $\Gamma(\emptyset)$ is determined by the set $J:=\{j\mid \varepsilon_j=-1\}$. More precisely, if $J=\{j_1<j_2<\ldots<j_s\}$, then $\Gamma=\phi_{j_1}\ldots\phi_{j_s}(\Gamma(\emptyset))$. We call the elements of $J$ the {\em folding positions} of $\Gamma$, and write $\Gamma=\Gamma(J)$. 

Throughout this paper, we will use $J$ and $\Gamma=\Gamma(J)$ interchangeably. For instance, according to the above discussion, we have
\[\phi_i(\Gamma(J))=\Gamma(J\triangle\{i\})\,,\]
where $\triangle$ denotes the symmetric difference of sets. Hence, it makes sense to define the folding operator $\phi_i$ on $J$ (compatibly with the action of $\phi_i$ on $\Gamma(J)$) by $\phi_i\::\:J\mapsto J\triangle\{i\}$. 

Using the notation above and the fact that $\gamma_{j_i}=r_{j_1}\ldots r_{j_{i-1}}(\beta_{j_i})$, we have 
\[{t}_{j_1}={r}_{j_1},\;{t}_{j_2}={r}_{j_1}{r}_{j_2}{r}_{j_1}, \; {t}_{j_3}={r}_{j_1}{r}_{j_2}{r}_{j_3}{r}_{j_2}{r}_{j_1},\:\ldots\,;\]
recall that $r_i=s_{\beta_i}$. In particular, we have
\begin{equation}\label{reflection}{r}_{j_1}\ldots {r}_{j_i}=({t}_{j_1}\ldots {t}_{j_i})^{-1}\,,\end{equation}
for $i=1,\ldots,s$. 

\begin{remark}\label{finiteness}  Although a folding $\Gamma$ of $\Gamma(\emptyset)$ is an infinite sequence if the root system is infinite, we are, in fact, always working with finite objects. Indeed, we are examining $\Gamma$ by considering only one root at a time.
\end{remark}

Given a folding $\Gamma$ of $\Gamma(\emptyset)$, we associate to each pair in $\Gamma$ (or the corresponding index $i$ in $I$) an integer $l_i$, which we call {\em level}; the sequence $L=L(\Gamma):=\{l_i\}_{i\in I}$ will be called the {\em level sequence} of $\Gamma$. The definition is as follows:
\begin{equation}\label{levelseq}l_i:=\delta+\stacksum{j<i,\varepsilon_j=1}{\gamma_j=\pm\gamma_i} {\rm sgn}(\gamma_j)\,,\end{equation}
where
\[\delta:=\casetwo{0}{\gamma_i>0}{-1}\]
We make the convention that the sum is 0 if it contains no terms. The definition makes sense since the sum is always finite. In particular, we have the level sequence $L_\emptyset=L(\Gamma(\emptyset)):=\{l_i^\emptyset\}_{i\in I}$ of $\Gamma(\emptyset)$. 

We will now consider certain affine reflections corresponding to foldings $\Gamma$ of $\Gamma(\emptyset)$. Let $\widehat{t}_i:=s_{|\gamma_i|,l_i}$; recall that the latter is the reflection in the affine hyperplane $H_{|\gamma_i|,l_i}$. In particular, we have the affine reflections $\widehat{r}_i:=s_{\beta_i,l_i^\emptyset}$ corresponding to $\Gamma(\emptyset)$.

\begin{definition}\label{defweighted} Given $J=\{j_1<j_2<\ldots<j_s\}\subseteq I$ and $\Gamma=\Gamma(J)$, we let 
\[\mu=\mu(\Gamma)=\mu(J):=\widehat{r}_{j_1}\ldots \widehat{r}_{j_s}(\lambda)\,,\]
and call $\mu$ the {\em weight} of $\Gamma$ (respectively $J$). We also let $\kappa(J)=\kappa(\Gamma):=r_{j_1}\ldots r_{j_s}$ (recall that $r_i:=s_{\beta_i}$).
\end{definition}

Given a root $\alpha$, we will use the following notation:
\begin{align}&I_\alpha=I_\alpha(\Gamma):=\{i\in I\mid  \gamma_i=\pm\alpha\}\,,\;\;\;\,L_\alpha=L_\alpha(\Gamma):=\{l_i\mid  i\in I_\alpha\}\,,\\
&\widehat{I}_\alpha=\widehat{I}_\alpha(\Gamma):=I_{\alpha}\cup\{\infty\}\,,\;\;\;\,\widehat{L}_\alpha=\widehat{L}_\alpha(\Gamma):=L_\alpha\cup\{l_\alpha^\infty\}\,,\nonumber\end{align}
where $l_\alpha^\infty:=\langle\mu(\Gamma),\alpha^\vee\rangle$. 

The following proposition is our main technical result, which relies heavily on the defining properties of $\lambda$-chains.

\begin{proposition}\label{hyperplane} Let $\Gamma=\Gamma(J)$ for some $J=\{j_1<j_2<\ldots<j_s\}\subseteq I$, and let $j_p<j\le j_{p+1}$ (the first or the second inequality is dropped if $p=0$ or $p=s$, respectively). Using the notation above, we have 
\[H_{|\gamma_j|,l_j}=\widehat{r}_{j_1}\ldots \widehat{r}_{j_p}(H_{\beta_j,l_j^\emptyset})\,.\]
\end{proposition}

\begin{proof}
It suffices to consider $p=s$. Let us define the roots $\delta_0:=\beta_j$ and
\[\delta_l:=r_{j_{s-l+1}}r_{j_{s-l+2}}\ldots r_{j_s}(\beta_j)\,,\]
for $l=1,\ldots,s$; note that $\delta_s=\gamma_j$. Throughout this proof, we use the maps $N_i$ and $N_{ij}$ defined at the beginning of Section \ref{sec:lambdachains} to count roots in $\Gamma(\emptyset)$. We need to show that
\[\widehat{r}_{j_1}\ldots \widehat{r}_{j_s}(H_{\beta_j,l_j^\emptyset})=H_{\delta_s,k}\,,\]
where, according to the definitions of folding operators and of levels, the integer $k:={\rm sgn}(\delta_s)l_j$ can be expressed as
\[k={\rm sgn}(\delta_s)N_{j_1}(|\delta_s|)+{\rm sgn}(\delta_{s-1})N_{j_1j_2}(|\delta_{s-1}|)+\ldots+N_{j_sj}(\delta_0)+\frac{1-{\rm sgn}(\delta_s)}{2}\,.\]
Indeed, by examining the segment of $\Gamma(\emptyset)$ between $j_l$ and $j_{l+1}$ for $l=1,\ldots,s$ (we let $j_{s+1}:=j$), we note that only the root $|\delta_{s-i}|$ gets changed to $\pm\delta_s$ in $\Gamma$; more precisely, the sign of the corresponding root in $\Gamma$ is ${\rm sgn}(\delta_s){\rm sgn}(\delta_{s-i})$. The definition of levels now leads us to the above formula for $k$. 

We now use induction on $s$, which starts at $s=0$. Let us assume that the statement to prove holds for $s-1$. Based on the discussion above, this means that
\[\widehat{r}_{j_2}\ldots \widehat{r}_{j_s}(H_{\beta_j,l_j^\emptyset})=H_{\delta_{s-1},k'}\,,\]
where
\[k'={\rm sgn}(\delta_{s-1})N_{j_2}(|\delta_{s-1}|)+{\rm sgn}(\delta_{s-2})N_{j_2j_3}(|\delta_{s-2}|)+\ldots+N_{j_sj}(\delta_0)+\frac{1-{\rm sgn}(\delta_{s-1})}{2}\,.\]
An easy linear algebra computation shows that
\[\widehat{r}_{j_1}(H_{\delta_{s-1},k'})=H_{\delta_s,k''}\,,\]
where
\[k''=k'-l_{j_1}^\emptyset\langle\beta_{j_1},\delta_{s-1}^\vee\rangle=k'-N_{j_1}(\beta_{j_1})\langle\beta_{j_1},\delta_{s-1}^\vee\rangle\,.\]
Let us now substitute into this formula the expression for $k'$ given by the induction hypothesis, and compare the result with the expression for $k$ given above. 

Let us first assume that $\delta_{s-1}\ne\pm\beta_{j_1}$. It turns out that verifying $k=k''$ amounts to proving the following identity:
\[{\rm sgn}(\delta_{s})N_{j_1}(|\delta_s|)-\frac{{\rm sgn}(\delta_{s})}{2}={\rm sgn}(\delta_{s-1})N_{j_1}(|\delta_{s-1}|)-\frac{{\rm sgn}(\delta_{s-1})}{2}-N_{j_1}(\beta_{j_1})\langle\beta_{j_1},\delta_{s-1}^\vee\rangle\,.\]
Let us now multiply both sides by ${\rm sgn}(\delta_{s-1})$, then substitute $\delta_s$ with $r_{j_1}(\delta_{s-1})$, and finally write $\alpha$ for $|\delta_{s-1}|$, $\beta$ for $\beta_{j_1}$, and $i$ for $j_1$. The identity to prove becomes precisely the one in Lemma~\ref{countroots}.

In the special case $\delta_{s-1}=\pm\beta_{j_1}$, we need to correct the identity above by adding ${\rm sgn}(\delta_{s-1})$ to its right-hand side. This is proved in a similar way, based on the correction of Lemma~\ref{countroots} in the case $\alpha=\beta$, which amounts to adding $1$ to the right-hand side of the corresponding formula. 
\end{proof}

\begin{corollary}\label{affinerefl}
Let $\Gamma=\Gamma(J)$ for some $J=\{j_1<j_2<\ldots<j_s\}\subseteq I$. Using the notation above, we have 
\[\widehat{t}_{j_1}=\widehat{r}_{j_1},\;\widehat{t}_{j_2}=\widehat{r}_{j_1}\widehat{r}_{j_2}\widehat{r}_{j_1}, \; \widehat{t}_{j_3}=\widehat{r}_{j_1}\widehat{r}_{j_2}\widehat{r}_{j_3}\widehat{r}_{j_2}\widehat{r}_{j_1},\:\ldots\]
In particular, we have
\[\widehat{r}_{j_1}\ldots \widehat{r}_{j_i}=(\widehat{t}_{j_1}\ldots \widehat{t}_{j_i})^{-1}\,,\]
for $i=1,\ldots,s$. Furthermore, if $\Gamma'=\phi_i(\Gamma)$, then $\mu(\Gamma')=\widehat{t}_i(\mu(\Gamma))$.
\end{corollary}

\begin{proof}
The first part follows from Proposition~\ref{hyperplane} by applying the following basic result: if $H_{\beta,m}=\widehat{r}_1\ldots\widehat{r}_q(H_{\alpha,k})$ for some affine reflections $\widehat{r}_1,\ldots, \widehat{r}_q$, then $s_{\beta,m}=\widehat{r}_1\ldots \widehat{r}_q s_{\alpha,k}\widehat{r}_q\ldots \widehat{r}_1$. The rest of the corollary follows easily from the first part.
\end{proof}

The next proposition shows that all inner products of $\mu(\Gamma)$ with positive roots can be easily read off from the level sequence $L(\Gamma)=(l_i)_{i\in I}$. Recall that, according to (\ref{levelseq}), the latter is computed by applying a simple counting procedure to the sequence of pairs in $\Gamma$. 

\begin{proposition}\label{lastlevel} Given a positive root $\alpha$, let $m:=\max\,I_\alpha(\Gamma)$, assuming that $I_\alpha(\Gamma)\ne\emptyset$. Then we have
\[\langle \mu(\Gamma),\alpha^\vee \rangle=\casethreeother{l_m+1}{\varepsilon_m\gamma_m>0\textrm{ and } t_{j_1}\ldots t_{j_s}(\alpha)>0}{l_m-1}{\varepsilon_m\gamma_m<0\textrm{ and } t_{j_1}\ldots t_{j_s}(\alpha)<0}{l_m}\]
On the other hand, if $I_\alpha(\Gamma)=\emptyset$, then we have
\[\langle \mu(\Gamma),\alpha^\vee \rangle=\casetwoex{0}{t_{j_1}\ldots t_{j_s}(\alpha)>0}{-1}{t_{j_1}\ldots t_{j_s}(\alpha)<0}\]
\end{proposition}

\begin{proof} We prove this result by induction on $s$, which starts at $s=0$. If $s>0$, we assume that the result holds for $\Gamma':=\Gamma(J\setminus\{i\})$, where $i:=j_1$. Let 
\[\beta:=\beta_i\,,\;\;\sigma:=\mathrm{sgn}(s_{\beta}(\alpha))\,,\;\;\mu:=\mu(\Gamma)\,,\;\;\mu':=\mu(\Gamma')\,,\;\;L'=(l_j')_{j\in I}:=L(\Gamma')\,.\]
We have
\begin{align}\label{weightcomp}\langle\mu,\alpha^\vee\rangle&=\langle\widehat{r}_i(\mu'),\alpha^\vee\rangle=\langle s_{\beta}(\mu'),\alpha^\vee\rangle+l_i^\emptyset\langle\beta,\alpha^\vee\rangle \\&=\sigma\,\langle\mu',|s_{\beta}(\alpha)|^\vee\rangle+l_i^\emptyset\langle\beta,\alpha^\vee\rangle \,.\nonumber\end{align}
Recall that $t_{j_1}\ldots t_{j_s}=r_{j_s}\ldots r_{j_1}$, by (\ref{reflection}). We have the following cases.

\underline{Case 1.1: $I_\alpha(\Gamma)\ne\emptyset$ and $m\ge i$.} We have $m=\max\,I_\alpha(\Gamma)=\max\,I_{s_\beta(\alpha)}(\Gamma')$. Based on (\ref{weightcomp}) and induction, we have
\[\langle\mu,\alpha^\vee\rangle=\sigma\,(l_m'+\delta')+l_i^\emptyset\langle\beta,\alpha^\vee\rangle \,,\]
where $\delta'\in\{-1,0,1\}$ is the correction term in the formula for $\langle\mu',|s_{\beta}(\alpha)|^\vee\rangle$. On the other hand, an easy linear algebra computation shows that $\widehat{r}_i$ maps the affine hyperplane $H_{|s_\beta(\alpha)|,l_m'}$ to 
\[H_{\sigma\alpha,l_m'-l_i^\emptyset\langle\beta,|s_\beta(\alpha)|^\vee\rangle}=H_{\alpha,\sigma l_m'+l_i^\emptyset\langle\beta,\alpha^\vee\rangle}\,.\]
 But then, by Proposition~\ref{hyperplane}, we have $l_m=\sigma l_m'+l_i^\emptyset\langle\beta,\alpha^\vee\rangle$. Hence, all we have to prove is that $\sigma \delta'=\delta$, where $\delta$ is the correction term in the formula for $\langle\mu,\alpha^\vee\rangle$. This follows from the fact that, if $(\gamma,\varepsilon)$ is the pair indexed by $m$ in $\Gamma'$, then we have
\[\mathrm{sgn}(s_\beta(\gamma))=\sigma\,\mathrm{sgn}(\gamma)\,,\;\;\;\mathrm{sgn}(r_{j_s}\ldots r_{j_1}(\alpha))=\sigma\,\mathrm{sgn}(r_{j_s}\ldots r_{j_2}(|s_\beta(\alpha)|))\,.\]

\underline{Case 1.2: $I_\alpha(\Gamma)\ne\emptyset$ and $m<i$.} It is clear that $\alpha\ne\beta$. We have two subcases, corresponding to $I_{s_\beta(\alpha)}(\Gamma')$ nonempty and empty. In the first subcase, letting $m':=\max\,I_{s_\beta(\alpha)}(\Gamma')$, we have $m'<i$, and, thus, $l_{m'}=N_i(|s_\beta(\alpha)|)-1$; here the counting function $N_i$ refers to our fixed $\lambda$-chain. Based on (\ref{weightcomp}), induction, and Lemma~\ref{countroots}, we have
\begin{align}\label{mucount}\langle\mu,\alpha^\vee\rangle&=\sigma\,(N_i(|s_\beta(\alpha)|)-1+\delta')+\langle\beta,\alpha^\vee\rangle N_i(\beta)\\&=N_i(\alpha)-\frac{\sigma+1}{2}+\sigma \delta'\,,\nonumber\end{align}
where $\delta'$ is as above. On the other hand, note that $l_m=N_i(\alpha)-1$, so we have to prove
\begin{equation}\label{toprove}
\langle\mu,\alpha^\vee\rangle=N_i(\alpha)-1+\delta\,.
\end{equation}
Hence, it remains to show that $(1-\sigma)/2+\sigma \delta'=\delta$, where $\delta$ is as above. This follows from the fact that 
\begin{equation}\label{defeps}\delta=\casetwoc{1}{r_{j_s}\ldots r_{j_1}(\alpha)>0}{0}\;\delta'=\casetwo{1}{\sigma\,r_{j_s}\ldots r_{j_1}(\alpha)>0}{0}\end{equation}
The subcase $I_{s_\beta(\alpha)}(\Gamma')=\emptyset$ is now immediate, since (\ref{mucount}) still holds, with $\delta'$ defined as in (\ref{defeps}); indeed, in this case the induction is based on the second formula in the proposition to be proved. 

\underline{Case 2: $I_\alpha(\Gamma)=\emptyset$.} This case is easily reduced to Case 1.2 above, since the formula to be proved can still be written as in (\ref{toprove}), with $\delta$ defined as in (\ref{defeps}). 
\end{proof}

\begin{remarks}\label{graphrep} (1) Since $\gamma_\infty=r_{j_1}\ldots r_{j_s}(\rho)$, we have ${\rm sgn}(t_{j_1}\ldots t_{j_s}(\alpha))={\rm sgn}(\langle\gamma_\infty,\alpha^\vee\rangle)$. Hence, Proposition~\ref{lastlevel} can be restated in terms of $\langle\gamma_\infty,\alpha^\vee\rangle$. 

(2) It is often useful to use the following graphical representation. Let $\widehat{I}_\alpha=\{i_1<i_2<\ldots<i_n=m<i_{n+1}=\infty\}$, and let us define the continuous piecewise-linear function $g_\alpha\::\:[0,n+\frac{1}{2}]\rightarrow{\mathbb R}$ by
\[g_\alpha(0)=-\frac{1}{2}\,,\;\;g_\alpha'(x)=\casethree{\mathrm{sgn}(\gamma_{i_k})}{x\in(k-1,k-\frac{1}{2}),\;\:k=1,\ldots,n}{\varepsilon_{i_k}\mathrm{sgn}(\gamma_{i_k})}{x\in(k-\frac{1}{2},k),\;\:k=1,\ldots,n}{{\rm sgn}(\langle\gamma_\infty,\alpha^\vee\rangle)}{x\in(n,n+\frac{1}{2})}\]
Then, according to the defining relation (\ref{levelseq}) and  Proposition~\ref{lastlevel}, we have 
\[l_{i_k}=g_\alpha\left(k-\frac{1}{2}\right),\;k=1,\ldots,n\,,\;\;\; \mbox{and}\;\;\; l_\alpha^\infty:=\langle\mu(\Gamma),\alpha^\vee\rangle=g_\alpha\left(n+\frac{1}{2}\right)\,.\]
For instance, assume that the entries of $\Gamma$ indexed by the elements of $I_\alpha$ are $(\alpha,-1)$, $(-\alpha,1)$, $(\alpha,1)$, $(\alpha,1)$, $(\alpha,-1)$, $(-\alpha,1)$, $(\alpha,-1)$, $(\alpha,1)$, in this order; also assume that ${\rm sgn}(\langle\gamma_\infty,\alpha^\vee\rangle)=1$. The graph of $g_\alpha$ is shown in Figure~\ref{fig:walls}.
\begin{figure}[ht]
\mbox{\epsfig{file=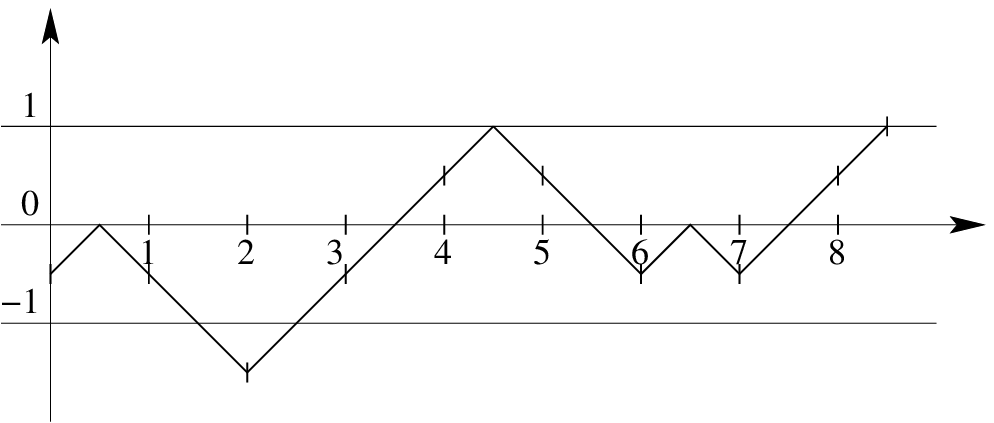}}
\caption{}
\label{fig:walls}
\end{figure}
\end{remarks}

\section{Admissible Subsets and Admissible Foldings}\label{sec:folded}

We will now define some special foldings $\Gamma(J)$ of $\Gamma(\emptyset)$. The notation is the same as in Section \ref{sec:chains-foldings}.

\begin{definition}\label{admsubset}
An {\em admissible subset} is a finite subset of $I$ (possibly empty), that is, $J=\{j_1<j_2<\ldots<j_s\}$, such that we have the following saturated chain in the Bruhat order on $W$:
\[1\lessdot r_{j_1}\lessdot r_{j_1}r_{j_2}\lessdot \ldots \lessdot r_{j_1}r_{j_2}\ldots r_{j_s}=\kappa(J)\,.\]
If $J$ is an admissible subset, we will call $\Gamma=\Gamma(J)$ an {\em admissible folding} (of $\Gamma(\emptyset)$). We denote by ${\mathcal A}$ the collection of all admissible subsets corresponding to our fixed $\lambda$-chain.
\end{definition}

\begin{remark}\label{lskey} The Weyl group element $\kappa(J)$ is a generalization of the {\em left-key} of semistandard Young tableaux, which was defined by Lascoux-Sch\"utzenberger \cite{lasksb}, and which is relevant to {\em Demazure characters}. Indeed, for finite types, a Demazure character formula in terms of $\kappa(J)$ appears in \cite{Le}.
\end{remark}

According to (\ref{reflection}), we have the following intrinsic criterion for $\Gamma(J)$ to be an admissible folding.

\begin{corollary}\label{intrinsic}
Given arbitrary $J=\{j_1<j_2<\ldots<j_s\}\subseteq I$, $\Gamma(J)$ is an admissible folding if and only if
\[1\lessdot t_{j_1}\lessdot t_{j_1}t_{j_2}\lessdot \ldots \lessdot t_{j_1}t_{j_2}\ldots t_{j_s}\,.\]
\end{corollary}

In this section, we will prove some results about admissible foldings. Let us fix an admissible subset $J=\{j_1<j_2<\ldots<j_s\}$, and let $\Gamma$ be the associated admissible folding; also, let $L(\Gamma)=(l_i)_{i\in I}$, $I_\alpha=I_\alpha(\Gamma)$, and $L_\alpha=L_\alpha(\Gamma)$. Recall the notation $L(\Gamma(\emptyset)):=\{l_i^\emptyset\}_{i\in I}$. We start with a basic result involving the Bruhat order.

\begin{lemma}\label{bruhat} Let $u \lessdot w = s_\beta u$ be a covering relation in the Bruhat order, where $\beta$ is some positive root, and $\alpha$ is a simple root. Assume that we have $u(\alpha)>0$ and $w(\alpha)<0$. Then $\beta = u(\alpha)$.
\end{lemma}

\begin{proof} For a simple root $\alpha$, the following conditions are equivalent: (1) $w(\alpha) < 0$; (2) $w s_\alpha < w$; (3) there exists a reduced decomposition for $w$ that ends with $s_\alpha$.  Let us pick such a reduced decomposition $w = s_{i_1} \ldots s_{i_l}$, so $\alpha = \alpha_{i_l}$. All elements $u$ that are covered by $w$ in the Bruhat order are obtained by skipping one term in the reduced decomposition.  We know that $u(\alpha)>0$, so the element $u$ cannot have a reduced decomposition that ends with $s_\alpha = s_{i_l}$. Thus $u$ is obtained from $w$ by skipping the last term $s_{i_l}$ in the reduced decomposition. So  $w = u s_\alpha = s_\beta u$, where $\beta=u(\alpha)$. 
\end{proof}

\begin{lemma}\label{existroot} Assume that $r_{j_a}\ldots r_{j_1}(\alpha)>0$ and $r_{j_b}\ldots r_{j_1}(\alpha)<0$ for some simple root $\alpha$ and $0\le a<b$ (if $a=0$, then the first condition is void). Then there exists $i$ with $a\le i<b$ such that $\gamma_{j_{i+1}}=\alpha$.
\end{lemma}

\begin{proof}
Find $i$ with $a\le i<b$ such that $r_{j_i}\ldots r_{j_1}(\alpha)>0$ and $r_{j_{i+1}}\ldots r_{j_1}(\alpha)<0$. By Lemma~\ref{bruhat}, we have $\beta_{j_{i+1}}=r_{j_i}\ldots r_{j_1}(\alpha)$. This means that $\gamma_{j_{i+1}}=r_{j_1}\ldots r_{j_i}(\beta_{j_{i+1}})=\alpha$.
\end{proof}

\begin{proposition}\label{folding} If $\varepsilon_i=-1$, then $\gamma_i>0$.
\end{proposition}

\begin{proof} We can assume that $i=j_s$. Since $\gamma_{i}=w(\beta_{i})$, where $w=r_{j_1}\ldots r_{j_{s-1}}$, we need to prove that the right-hand side is a positive root. This follows from the fact that $w< ws_{\beta_i}$, according to a well-known property of the Bruhat order on a Coxeter group, cf. \cite[Proposition~5.7]{Hum}. 
\end{proof}

\begin{proposition}\label{signchange} Assume that $\alpha$ is a simple root for which $I_\alpha\ne \emptyset$. Let $m\in I_\alpha$ be either the minimum of $I_\alpha$, or an element for which its predecessor $k$ (in $I_\alpha$) satisfies $(\gamma_k,\varepsilon_k)=(\alpha,1)$. Then we have $\gamma_m=\alpha$. 
\end{proposition}

\begin{proof}
Assume that $\gamma_m=-\alpha$. Recall that the corresponding admissible subset is $J=\{j_1<j_2<\ldots<j_s\}$. Note that $m\not\in J$, based on Proposition~\ref{folding} (indeed, we must have $\varepsilon_m=1$). Let us define the index $b$ by the condition $j_b<m<j_{b+1}$ (possibly, $b+1=s$, in which case the second inequality is dropped). We also define the index $a$ by setting $a:=0$ if $m=\min\,I_\alpha$, and by the condition $j_a<k<j_{a+1}$, otherwise (if $a=0$ in the second case, the corresponding first inequality is dropped). We clearly have $r_{j_1}\ldots r_{j_{b}}(\beta_m)=-\alpha$, which implies $r_{j_b}\ldots r_{j_1}(\alpha)<0$. If $a\ne 0$, we also have $r_{j_1}\ldots r_{j_{a}}(\beta_k)=\alpha$, so $r_{j_a}\ldots r_{j_1}(\alpha)>0$. But then Lemma~\ref{existroot} applies and leads to a contradiction.
\end{proof}

\begin{proposition}\label{lastneg} Assume that, for some simple root $\alpha$, we have either $I_\alpha=\emptyset$, or $(\gamma_m,\varepsilon_m)=(\alpha,1)$ for $m=\max\,I_\alpha$. Then we have $\langle\gamma_\infty,\alpha^\vee\rangle>0$.
\end{proposition}

\begin{proof} Assume that the conclusion fails, which means that $r_{j_s}\ldots r_{j_1}(\alpha)<0$ (cf.\ Remark \ref{graphrep} (2)). Define the index $a$ by setting $a:=0$ if $I_\alpha=\emptyset$, and by the condition $j_a<m<j_{a+1}$, otherwise (if $a=0$ or $a=s$ in the second case, one of the two inequalities is dropped). If $a\ne 0$, we have $r_{j_1}\ldots r_{j_a}(\beta_m)=\alpha$, so $r_{j_a}\ldots r_{j_1}(\alpha)>0$. It means that the hypotheses of Lemma~\ref{existroot} are satisfied for $b:=s$. This lemma now leads to a contradiction. 
\end{proof}

\begin{proposition}\label{maximum} Assume that $s>0$ and let $J':=J\setminus\{j_s\}$. We have
\[\mu(J')-\mu(J)=\left(\langle\lambda,\beta_{j_s}^\vee\rangle-l_{j_s}^\emptyset\right)\gamma_{j_s}\,.\]
In other words, $\mu(J')-\mu(J)$ is a positive multiple of a positive root. In particular, if $J\ne\emptyset$ then $\lambda-\mu(J)\ne 0$ belongs to the positive cone in the root lattice.
\end{proposition}

\begin{proof}
Using the definitions and notation in Section \ref{sec:chains-foldings}, we have
\begin{align*}&\mu(J')-\mu(J)=\widehat{r}_{j_1}\ldots\widehat{r}_{j_{s-1}}(\lambda)-\widehat{r}_{j_1}\ldots\widehat{r}_{j_{s-1}}(\widehat{r}_{j_{s}}(\lambda))\\&={r}_{j_1}\ldots{r}_{j_{s-1}}(\lambda-\widehat{r}_{j_{s}}(\lambda))=\left(\langle\lambda,\beta_{j_s}^\vee\rangle-l_{j_s}^\emptyset\right){r}_{j_1}\ldots{r}_{j_{s-1}}(\beta_{j_s})=\left(\langle\lambda,\beta_{j_s}^\vee\rangle-l_{j_s}^\emptyset\right)\gamma_{j_s}\,.\end{align*}
The fact that $\gamma_{j_s}$ is a positive root is the content of Proposition \ref{folding}.
\end{proof}

\begin{proposition}\label{extremalelement} Given $\Gamma=\phi_{j_1}\ldots \phi_{j_s}(\Gamma(\emptyset))$ and $\mu={\mu}(\Gamma)$, the following 
hold.
\begin{enumerate} 
\item If $s=0$, then $\max\,L_\alpha<\langle\mu,\alpha^\vee\rangle$ for all simple roots $\alpha$ with $I_\alpha\ne\emptyset$. 
\item If $s>0$, then there is a simple root $\alpha$ in $\Gamma$ such that $I_\alpha\ne\emptyset$ and $\max\,L_\alpha>\langle\mu,\alpha^\vee\rangle$.
\end{enumerate}
\end{proposition}

\begin{proof} The first part follows directly from the definitions. If $s>0$, we can find a simple root $\alpha$ such that $r_{j_s}\ldots r_{j_1}(\alpha)<0$. Hence, Proposition~\ref{lastneg} applies, and, letting $m:=\max\,I_\alpha$, we have $\varepsilon_m\gamma_m=-\alpha$. But then $\langle\mu,\alpha^\vee\rangle=l_m-1$, by Proposition~\ref{lastlevel}, so the second statement is verified. 
\end{proof}

Let us now fix a simple root $\alpha$, and recall the description of the sets $\widehat{I}_\alpha(\Gamma)$ and $\widehat{L}_\alpha(\Gamma)$ based on the continuous piecewise-linear function $g_\alpha$, which was introduced in Remark \ref{graphrep} (2). We will rephrase some of the above results in a simple way in terms of $g_\alpha$, and we will deduce some consequences (to be used in the subsequent sections), which are formulated in the same language. Assume that $I_\alpha=\{i_1<i_2<\ldots<i_n\}$, so that $g_\alpha$ is defined on $[0,n+\frac{1}{2}]$, and let $M$ be the maximum of $g_\alpha$. Note first that the function $g_\alpha$ is determined by the sequence $(\sigma_1,\ldots,\sigma_{n+1})$, where $\sigma_j:=(\mathrm{sgn}(\gamma_{i_j}), \varepsilon_{i_j} \mathrm{sgn}(\gamma_{i_j}))$ for $1\le j\le n$, and $\sigma_{n+1}:= \mathrm{sgn}(\langle\gamma_{\infty},\alpha^\vee\rangle)$. We have the following restrictions.
\begin{enumerate}
\item[(C1)] $\sigma_j\in\{(1,1),\,(-1,-1),\,(1,-1)\}$ for $1\le j\le n$ (by Proposition~\ref{folding}).
\item[(C2)] $j=0$ or $\sigma_j=(1,1)$ implies $\sigma_{j+1}\in\{(1,1),\,(1,-1),\,1\}$ (by Propositions \ref{signchange} and \ref{lastneg}).
\end{enumerate}

\begin{corollary}\label{cor1} We have $M\ge 0$. If $g_\alpha(x)=M$, then $x=m+\frac{1}{2}$ for $0\le m\le n$, and $\sigma_{m+1}\in\{(1,-1),1\}$.
\end{corollary}

\begin{proof} The fact that $M\ge 0$ follows from Condition (C2) for $j=0$. Thus, $g_\alpha(0)=-\frac{1}{2}\ne M$. If $g_\alpha(m)=M$ for $m\in\{1,\ldots,n\}$, then $\sigma_m=(1,1)$ by Condition (C1). But then Condition (C2) leads to a contradiction. The last statement is obvious.
\end{proof}

\begin{corollary}\label{cor2} Assume that $M>0$, and let $m$ be such that $m+\frac{1}{2}=\min\,g_\alpha^{-1}(M)$. We have $m>0$, $\sigma_{m}=(1,1)$, and $g_\alpha(m-\frac{1}{2})=M-1$. Moreover, we have $g_\alpha(x)\le M-1$ for $0\le x\le m-\frac{1}{2}$.
\end{corollary}

\begin{proof}
If $m=0$, then $g_\alpha(m+\frac{1}{2})=0$, so $M=0$. It is also easy to see that both $\sigma_{m}=(-1,-1)$ and $\sigma_{m}=(1,-1)$ contradict the definition of $m$. Now assume that the last statement fails. Then there exists $1\le k\le m-1$ such that $g_\alpha(k-1)=M-\frac{1}{2}$ and $\sigma_k=(-1,-1)$. Condition (C2) implies $k\ge 2$ and $\sigma_{k-1}\in\{(-1,-1),\,(1,-1)\}$, which is a contradiction to the definition of $m$. 
\end{proof}

The next corollary can be proved in a similar way to Corollary \ref{cor2}. 

\begin{corollary}\label{cor3} Assume that $M>g_\alpha(n+\frac{1}{2})$, and let $k$ be such that $k-\frac{1}{2}=\max\,g_\alpha^{-1}(M)$. We have $k\le n$, $\sigma_{k+1}=(-1,-1)$, and $g_\alpha(k+\frac{1}{2})=M-1$. Moreover, we have $g_\alpha(x)\le M-1$ for $k+\frac{1}{2}\le x\le n+\frac{1}{2}$.
\end{corollary}

\section{Root Operators}\label{sec:rootoperators}

We will now define {\em root operators} on the collection ${\mathcal A}$ of admissible subsets corresponding to our fixed $\lambda$-chain. Let $J$ be such an admissible subset, let $\Gamma$ be the associated admissible folding, and $L(\Gamma)=(l_i)_{i\in I}$ its level sequence, denoted as in Section \ref{sec:chains-foldings}. Also recall from Section \ref{sec:chains-foldings} the definitions of the finite sequences $I_\alpha(\Gamma)$, $\widehat{I}_\alpha(\Gamma)$, $L_\alpha(\Gamma)$, and $\widehat{L}_\alpha(\Gamma)$, where $\alpha$ is a root. 

We will first define a partial operator $F_p$ on admissible subsets $J$ for each $p$ in $\{1,\ldots,r\}$, that is, for each simple root $\alpha_p$. Let $p$ in $\{1,\ldots,r\}$ be fixed throughout this section. Let $M=M(\Gamma)=M(\Gamma,p)=M(J,p)$ be the maximum of the finite set of integers $\widehat{L}_{\alpha_p}(\Gamma)$. We know that $M\ge 0$ from Corollary \ref{cor1}.  
 Assume that $M>0$. Let $m=m_F(\Gamma)=m_F(\Gamma,p)$ be the minimum index $i$ in $I_{\alpha_p}(\Gamma)$ for which we have $l_i=M$.
If no such index exists, then $M=\langle \mu(\Gamma),\alpha_p^\vee\rangle$; in this case, we let $m=m_F(\Gamma)=m_F(\Gamma,p):=\infty$. Now let $k=k_F(\Gamma)=k_F(\Gamma,p)$ be the predecessor of $m$ in $\widehat{I}_{\alpha_p}(\Gamma)$. 
 By Corollary \ref{cor2}, this always exists and we have $l_k=M-1\ge 0$. 


Let us now define
\begin{equation}\label{defroot}
F_p(J):=\phi_k \phi_m(J)\,,
\end{equation}
where $\phi_\infty$ is the identity map. Note that the folding of $\Gamma(\emptyset)$ associated to $F_p(J)$, which will be denoted by $F_p(\Gamma)=(\{(\delta_i,\zeta_i)\}_{i\in I},\delta_\infty)$, is defined by a similar formula. More precisely, we have
\[(\delta_i,\zeta_i)=\casefour{(\gamma_i,\varepsilon_i)}{i<k\textrm{ or } i>m}{(\gamma_i,-\varepsilon_i)}{i=k}{(s_p(\gamma_i),\varepsilon_i)}{k<i<m}{(s_p(\gamma_i),-\varepsilon_i)}{i=m}\]
and 
\[\delta_\infty=\casetwoex{\gamma_\infty}{m\ne\infty}{s_p(\gamma_\infty)}{m=\infty}\]
In other words, based on Proposition~\ref{propf} below, we can say that applying the root operator $F_p$ amounts to performing a ``folding'' in position $k$, and, if $m\ne\infty$, an ``unfolding'' in position $m$. 

\begin{proposition}\label{propf} Given the above setup, the following hold.
\begin{enumerate}
\item If $m\ne\infty$, then $\gamma_m=\alpha_p$ and $\varepsilon_m=-1$.
\item We have $\gamma_k=\alpha_p$ and $\varepsilon_k=1$. 
\item We have ${\mu}(F_p(J))={\mu}(J)-\alpha_p\,.$
\end{enumerate}
\end{proposition}

\begin{proof} 


Let $\mu=\mu(\Gamma)$. The first two statements follow immediately from Corollaries \ref{cor1} and \ref{cor2}. For the third statement, note that, by Corollary \ref{affinerefl}, the weight of $F_p(J)$ is $\widehat{t}_k\widehat{t}_m(\mu)$ if $m\ne\infty$, and $\widehat{t}_k(\mu)$ otherwise. Using the formula $\widehat{t}_k(\nu)= s_p(\nu)+l_k\alpha_p$, and the similar one for $\widehat{t}_m$, we compute (in both cases)
\[{\mu}(F_p(J))=\mu+(l_k-M)\alpha_p=\mu-\alpha_p\,.\]
\end{proof}

We now intend to define a partial inverse $E_p$ to $F_p$. Assume that $M>\langle \mu(\Gamma),\alpha_p^\vee\rangle$.  Let $k=k_E(\Gamma)=k_E(\Gamma,p)$ be the maximum index $i$ in $I_{\alpha_p}(\Gamma)$ for which we have $l_i=M$.
Note that such indices always exist, by Corollary \ref{cor3}. Now let $m=m_E(\Gamma)=m_E(\Gamma,p)$ be the successor of $k$ in $\widehat{I}_{\alpha_p}(\Gamma)$. Corollary \ref{cor3} implies that, if $m=\infty$, then we have $\langle \mu(\Gamma),\alpha_p^\vee \rangle=M-1$, while, otherwise, we have $l_m=M-1$. 
Finally, we define $E_p(J)$ by the same formula as $F_p(J)$, namely (\ref{defroot}). Hence, the folding of $\Gamma(\emptyset)$ associated to $E_p(J)$ is also defined in the same way as above. The following analog of Proposition~\ref{propf} is proved in a similar way, by invoking Corollaries \ref{cor1} and \ref{cor3}.

\begin{proposition}\label{prope} Given the above setup, the following hold.
\begin{enumerate}
\item We have $\gamma_k=\alpha_p$ and $\varepsilon_k=-1$.
\item If $m\ne\infty$, then $\gamma_m=-\alpha_p$ and $\varepsilon_m=1$. 
\item We have ${\mu}(E_p(J))={\mu}(J)+\alpha_p\,.$
\end{enumerate}
\end{proposition}




\begin{proposition}\label{closed} If $F_p(J)$ is defined, then it is also an admissible subset. Similarly for $E_p(J)$. Moreover, if $F_p(J)$ and $E_p(J)$ are both defined, then $\kappa(F_p(J))=\kappa(J)$; otherwise, we have $\kappa(F_p(J))=s_p \kappa(J)$. 
\end{proposition}

\begin{proof}
We consider $F_p$ first. The cases corresponding to $m\ne\infty$ and $m=\infty$ can be proved in similar ways, so we only consider the first case. Let $J=\{j_1<j_2<\ldots<j_s\}$, as usual, and, based on Proposition~\ref{propf} (1)-(2), let $a<b$ be such that 
\[j_{a}<k<j_{a+1}<\ldots<j_b=m<j_{b+1}\,;\]
if $a=0$ or $b+1>s$, then the corresponding indices $j_a$, respectively $j_{b+1}$, are missing. If $a+1=b$, there is nothing to prove, so we assume $a+1<b$. 

We use the criterion in Corollary \ref{intrinsic}. Based on it, we have
\begin{equation}\label{chain1}
1\lessdot t_{j_1}\lessdot t_{j_1}t_{j_2}\lessdot \ldots \lessdot t_{j_1}t_{j_2}\ldots t_{j_s}\,.\end{equation}
Let $w:=t_{j_1}\ldots t_{j_a}$; if $a=0$, then set $w:=1$. We need to prove that
\begin{equation}\label{chain2}
w\lessdot ws_p\lessdot ws_pt_{j_{a+1}}'\lessdot\ldots\lessdot ws_pt_{j_{a+1}}'\ldots t_{j_{b-1}}'=t_{j_1}\ldots t_{j_b}\,,\end{equation}
where $t_{j_i}':=s_pt_{j_i}s_p$ for $i=a+1,\ldots,b-1$. Indeed, the pair with index $i$ in $F_p(\Gamma)$ is $(s_p(\gamma_i),\varepsilon_i))$, for $k<i<m$. Note that (\ref{chain2}) also implies that $\kappa(F_p(J))=\kappa(J)$. 

Recall that $r_{j_1}\ldots r_{j_i}=t_{j_i}\ldots t_{j_1}$, by (\ref{reflection}). On the other hand, based on Proposition~\ref{propf} (2), we have $r_{j_1}\ldots r_{j_a}(\beta_k)=\alpha_p$, which implies $w(\alpha_p)>0$. Hence, we have $w\lessdot ws_p$, which gives us the first covering relation in (\ref{chain2}). But we also have $w\lessdot wt_{j_{a+1}}$, by (\ref{chain1}), and $t_{j_{a+1}}\ne s_p$, by the choice of $k$. Based on a well-known property of the Bruhat order on a Coxeter group \cite[Theorem~1.1~(IV)~(iii)]{Deo1}, we deduce that $wt_{j_{a+1}}\lessdot wt_{j_{a+1}}s_p$ and $ws_p\lessdot wt_{a_{j+1}}s_p=ws_p t_{a_{j+1}}'$. The latter statement gives us the second covering relation in (\ref{chain2}). We can proceed in this way until we get the whole chain (\ref{chain2}). The proof of the result for $E_p(J)$ is completely similar, based on Proposition~\ref{prope}.
\end{proof}

\begin{proposition}\label{bijection} {\rm (1)} Assume that $F_p(J)$ is defined. Then we have
\begin{align*}&M(F_p(\Gamma))=M(\Gamma)-1>\langle{\mu}(F_p(\Gamma)),\alpha_p^\vee\rangle,\;\;\mathrm{and}\\
&k_E(F_p(\Gamma))=k_F(\Gamma),\;\;m_E(F_p(\Gamma))=m_F(\Gamma)\,.\end{align*}
Hence, $E_p(F_p(J))$ is defined and equal to $J$. 

{\rm (2)} Assume that $E_p(J)$ is defined. Then we have
\begin{align*}&M(E_p(\Gamma))=M(\Gamma)+1>0,\;\;\mathrm{and}\\
&k_F(E_p(\Gamma))=k_E(\Gamma),\;\;m_F(E_p(\Gamma))=m_E(\Gamma)\,.\end{align*}
Hence $F_p(E_p(J))$ is defined and equal to $J$. 

{\rm (3)} Let $a$ be maximal such that $F_p^a(J)$ is defined, and let $b$ be maximal such that $E_p^b(J)$ is defined. Then $a-b=\langle\mu,\alpha_p^\vee\rangle$.
\end{proposition}

The proof is straightforward based on the results already proved in this section, and is left to the interested reader. We note that that it is convenient to use the graphical representation described in Remark \ref{graphrep} (2).

\section{Admissible Subsets Form a Semiperfect Crystal}
\label{sec:mainthm}

In this section, we derive our main result. We start with the following definitions and lemma. Given an admissible folding $\Gamma=\Gamma(J)$, denoted as in Section \ref{sec:chains-foldings}, we denote by $\Gamma|_{\ge i}$ the sequence indexed by $j\ge i$ given by $i\mapsto\varepsilon_i\gamma_i$ and $j\mapsto(\gamma_j,\varepsilon_j)$ for $j>i$. Let us define
\[\varepsilon(J,p)=\varepsilon(\Gamma,p):=M(J,p)\,,\;\;\;\:\delta(J,p)=\delta(\Gamma,p):=\langle\mu(J),\alpha_p^\vee\rangle-M(J,p)\,,\]
where $M(J,p)$ was defined at the beginning of Section \ref{sec:rootoperators}.

\begin{lemma}\label{lemtiming} 
Let $\Delta$ and $\Gamma$ be two admissible foldings, with $L(\Delta)=(l_j)_{j\in I}$ and $L(\Gamma)=(l_j')_{j\in I}$. Assume that $\kappa(\Delta)=\kappa(\Gamma)$, that $E_p(\Delta)$ is defined, and that $\Delta_{\ge i}=\Gamma|_{\ge i}$ for $i=k_E(\Delta,p)$. Then $\delta(\Gamma,p)\le\delta(\Delta,p)$. Moreover, equality holds if and only if $l_i'=\max\,L_{\alpha_p}(\Gamma)$. 
\end{lemma}

\begin{proof}
Let $j_1,\ldots,j_s$ and $j_1',\ldots,j_t'$ be the folding positions of $\Delta$ and $\Gamma$, respectively, where
\begin{align*}
&j_1<\ldots<j_a= i<j_{a+1}<\ldots<j_s\,,\\
&j_1'<\ldots<j_b'\le i<j_{b+1}'=j_{a+1}<\ldots<j_t'=j_s\,.\end{align*}
Note that the pair in $\Delta$ indexed by $i$ is $(\alpha_p,-1)$, by Proposition~\ref{prope} (1). Hence, we have $r_{j_1}\ldots r_{j_a}(\beta_i)=r_{j_1'}\ldots r_{j_b'}(\beta_i)=-\alpha_p$, so $r_{j_s}\ldots r_{j_1}(\alpha_p)=r_{j_t'}\ldots r_{j_1'}(\alpha_p)$. By Proposition~\ref{lastlevel}, we have
\begin{align*}\delta(\Delta,p)&=\langle \mu(\Delta),\alpha_p^\vee\rangle-M(\Delta,p)=\langle \mu(\Delta),\alpha_p^\vee\rangle-l_i\\
&=\langle \mu(\Gamma),\alpha_p^\vee\rangle-l_i'\ge \langle \mu(\Gamma),\alpha_p^\vee\rangle-M(\Gamma,p)=\delta(\Gamma,p)\,.\end{align*}
\end{proof}

Recall that ${\mathcal A}$ is the collection of all admissible subsets corresponding to our fixed $\lambda$-chain.

\begin{theorem} 
The collection ${\mathcal A}$ of admissible subsets together with the root operators
form a semiperfect crystal.  Thus we have the following character formula:
\[\chi(\lambda)=\sum_{J\in{\mathcal A}}e^{\mu(J)}\,.\]
\end{theorem}

\begin{proof} As usual, we denote by $J$ a generic element of $\mathcal A$, by $\mu$ its weight, and by $\Gamma=\Gamma(J)$ the corresponding admissible folding. Also, $\alpha_p$ will be a generic simple root.

Axiom (A1) follows from Corollary \ref{cor1}. Axiom (A2) is the content of Proposition~\ref{closed} and Proposition~\ref{bijection} (1)-(2). Axiom (A3) is the content of Proposition~\ref{propf} (3) and Proposition~\ref{bijection} (1). 

According to Proposition~\ref{extremalelement},  $\emptyset$ is the unique admissible subset which is maximal with respect to all partial orders $\preceq_p$, for $p=1,\ldots,r$. Since the height of $\lambda=\mu(\emptyset)$ is strictly larger than the height of $\mu(J)$ for any other admissible subset $J$ (cf. Proposition \ref{maximum}), and since the height of $\mu(E_p(J))$ is larger by 1 than the height of $\mu(J)$, we conclude that $\emptyset$ is the maximum object of ${\mathcal A}$. Hence, Axiom (A5) is satisfied. 

The definition of a coherent timing pattern and the related verification of Axiom (A4) are analogous to those for LS chains in \cite[Theorem~8.3]{St}; nevertheless, there are some features specifically related to our setup, such as the reversal of the total order on the set in which the timing pattern takes values, and the use of Proposition~\ref{lastlevel} in Lemma~\ref{lemtiming} (used below). Assume that $\delta(J,p)<0$, so that $E_p(J)$ is defined. We define $t(J,p)=t(\Gamma,p):=k_E(\Gamma,p)$, where $k_E(\Gamma,p)$ was defined in Section \ref{sec:rootoperators} in connection with the root operator $E_p$. Assuming that $F_p(J)$ is defined, and applying Proposition~\ref{bijection} (1), we have
\[t(\Gamma,p)=k_E(\Gamma,p)\ge m_F(\Gamma,p)>k_F(\Gamma,p)=k_E(F_p(\Gamma),p)=t(F_p(\Gamma),p)\,.\]
By iteration, it follows that $\Gamma'|_{\ge t(\Gamma,p)}=\Gamma|_{\ge t(\Gamma,p)}$ for $\Gamma'\preceq_p\Gamma$. Therefore, given $\Delta\succeq_q\Gamma$ for $q\ne p$ with $\delta(\Delta,q)=\delta<0$ and $t=t(\Delta,q)\ge t(\Gamma,p)$, we have
\[\Delta|_{\ge t}=\Gamma|_{\ge t}=F_p(\Gamma)|_{\ge t}\,.\]
We have $\kappa(\Delta)=\kappa(\Gamma)=\kappa(F_p(\Gamma))$, by Proposition~\ref{closed}. Hence, by Lemma~\ref{lemtiming}, we have $\delta(F_p(\Gamma),q)\le \delta$, so there exists $\Delta'\succeq_q F_p(\Gamma)$ such that $\delta(\Delta',q)=\delta$. We claim that $t(\Delta',q)=t(\Delta,q)$. Let $t':=t(\Delta',q)$ and $t_*:=\max(t,t')$. We have 
\[\Delta'|_{\ge t_*}=F_p(\Gamma)|_{\ge t_*}=\Gamma|_{\ge t_*}=\Delta|_{\ge t_*}\,.\]
Also note that $\kappa(\Delta)=\kappa(F_p(\Gamma))=\kappa(\Delta')$, by Proposition~\ref{closed}. Assume that $t<t'$. Then the fact that $\delta(\Delta,q)=\delta(\Delta',q)$ implies, based on Lemma~\ref{lemtiming}, that the maximum of $L_{\alpha_q}(\Delta)$ is attained at $t'$ as well; but this contradicts the definition of $t$. The case $t'<t$ is similar. On the other hand, similar reasoning proves conversely that given $\Delta'\succeq_q F_p(\Gamma)$ such that $\delta(\Delta',q)=\delta<0$ and $t(\Delta',q)\ge t(\Gamma,p)$, there is $\Delta\succeq_q \Gamma$ such that $\delta(\Delta,q)=\delta$ and $t(\Delta,q)=t(\Delta',q)$.

The formula for $\chi(\lambda)$ now follows from Theorem~\ref{thmst}.
\end{proof}

\begin{corollary} (Littlewood-Richardson rule). We have
\[\chi(\lambda)\cdot \chi(\nu)=\sum\chi(\nu+\mu(J))\,,\]
where the summation is over all $J$ in ${\mathcal A}$ satisfying $\langle\nu+\mu(J),\alpha_p^\vee\rangle\ge M(J,p)$ for all $p=1,\ldots,r$.
\end{corollary}

\begin{proof} This follows immediately from Theorem~\ref{thmst}. There, the condition for $J$ to contribute to the summation was $\langle\nu,\alpha_p^\vee\rangle+\delta(J,p)\ge 0$, which is easily seen to be equivalent to the condition stated above.
\end{proof}

\begin{corollary} (Branching rule). Given $P\subseteq\{1,\ldots,r\}$, we have the following rule for decomposing $\chi(\lambda)$ as a sum of Weyl characters relative to $\Phi_P$:
\[\chi(\lambda)=\sum\chi(\mu(J);P)\,,\]
where the summation is over all $J$ in ${\mathcal A}$ satisfying $\langle \mu(J),\alpha_p^\vee\rangle=M(J,p)$ for all $p\in P$.
\end{corollary}

\begin{proof} Immediate, based on Corollary \ref{branching}.
\end{proof}

\section{Lakshmibai-Seshadri Chains}\label{sec:ls}

In this section, we explain the connection between our model and LS chains. We start with the relevant definitions, by closely following \cite[Section 8]{St}.

The Bruhat order on the orbit $W\lambda$ of a dominant or antidominant weight is the transitive closure of the relations
\[s_\alpha(\mu)<\mu\;\;\;\;\;\mbox{if}\;\;\;\;\;\langle\mu,\alpha^\vee\rangle >0\;\;\;(\mu\in W\lambda,\,\alpha\in\Phi^+)\,.\]
The Bruhat orders on $W\lambda$ and $-W\lambda$ are dual isomorphic; in fact, $\mu<\nu$ if and only if $-\nu<-\mu$. As usual, we write $\nu\lessdot\mu$ to indicate that $\mu$ covers $\nu$; this happens only if $\nu=s_\alpha(\mu)$ for some $\alpha\in\Phi^+$, but not conversely. Given $\pm\lambda\in \Lambda^+$ and a fixed real number $b$, one defines the $b$-{\em Bruhat order} $<_b$ as the transitive closure of the relations
\[s_\alpha(\mu)<_b\mu\;\;\;\;\;\mbox{if}\;\;\;\;\;s_\alpha(\mu)\lessdot\mu\;\;\:\mbox{and}\;\;\: b\langle\mu,\alpha^\vee\rangle\in{\mathbb Z}\;\;\;(\mu\in W\lambda,\,\alpha\in\Phi^+)\,.\]
Thus, $\mu$ covers $\nu$ in $b$-Bruhat order if and only if $\mu$ covers $\nu$ in the usual Bruhat order and $b(\mu-\nu)$ is an integer multiple of a root.

\begin{definition} Given $\pm\lambda\in \Lambda^+$, we say that a pair consisting of a chain $\mu_0<\mu_1<\ldots<\mu_l$ in the $W$-orbit of $\lambda$ and an increasing sequence of rational numbers $0<b_1<\ldots<b_l<1$ is a {\em Lakshmibai-Seshadri chain} (LS chain) if
\[\mu_0<_{b_1}\mu_1<_{b_2}\ldots<_{b_l}\mu_l\,.\]
\end{definition}

Following \cite{St}, we identify an LS chain (denoted as above) with the map $\gamma\,:\,(0,1]\rightarrow W\lambda$ given by $\gamma(t):=\mu_k$ for $b_k<t\le b_{k+1}$, where $k=0,\ldots,l$ and $b_0:=0$, $b_{l+1}:=1$. Note that the piecewise-constant left-continuous maps that correspond to LS chains can be characterized by the property
\[\gamma(t)\le_t\gamma(t^+)\;\;\:\mbox{for}\;\;\:0<t<1\,,\]
where $\gamma(t^+)$ denotes the right-hand limit of $\gamma$ at $t$. To each LS chain $\gamma$, we associate the continuous piecewise-linear path $\pi\,:\,[0,1]\rightarrow\hR$ given by
\[\pi(t):=\int_0^t\gamma(s)\,ds\,.\]
In other words, we define
\[\pi(t):=(t-b_{k})\mu_{k}+\sum_{i=0}^{k-1}(b_{i+1}-b_{i})\mu_{i}\,,\]
for $b_{k}\le t\le b_{k+1}$. 

The root operators $F_p$ and $E_p$ on LS chains (for a simple root $\alpha_p$) were defined by Littelmann \cite{Li1,Li2} as follows. Let $m_p$ be the minimum of the function $h_p\,:\,[0,1]\rightarrow {\mathbb R}$ given by $t\mapsto \langle\pi(t),\alpha_p^\vee\rangle$. If $m_p> -1$, then $E_p$ is undefined. Otherwise, we let $t_1\in [0,1]$ be minimal such that $h_p(t_1)=m_p$, and we let $t_0\in [0,t_1]$ be maximal such that $h_p(t)\ge m_p+1$ for $t\in[0,t_0]$. We define
\begin{equation}\label{defeppaths}E_p(\gamma)(t):=\casetwo{s_p(\gamma(t))}{t_0<t\le t_1}{\gamma(t)}\end{equation}
The definition of $F_p$ is similar. More precisely, if $h_p(1)-m_p<1$, then $F_p$ is undefined. Otherwise, we let $t_0\in [0,1]$ be maximal such that $h_p(t_0)=m_p$, and we let $t_1\in [t_0,1]$ be minimal such that $h_p(t)\ge m_p+1$ for $t\in [t_1,1]$. Given  the latter values for $t_0$ and $t_1$, we define $F_p(\gamma)$ by the same formula as $E_p(\gamma)$. 

\begin{remark}\label{defrootops} The above definition of $E_p$ and $F_p$ applies to any continuous piecewise-linear path $\pi$ with $\pi(0)=0$, if we replace the map $\gamma$ above with the left-hand derivative of $\pi$ (which is a piecewise-constant left-continuous map defined on $(0,1]$). 
\end{remark}

With each LS chain $\gamma$ (denoted as above), one can associate the {\em dual} LS chain $\gamma^*$, which is defined as
\[-\mu_l<_{a_l}\ldots<_{a_2}-\mu_1<_{a_1}-\mu_0\,,\]
where $a_i:=1-b_i$. It is easy to see that $E_p(\gamma^*)=F_p(\gamma)^*$ and $F_p(\gamma^*)=E_p(\gamma)^*$. 

The paths corresponding to LS chains are special cases of {\em Littelmann paths}. The latter were defined by Littelmann \cite{Li2}, still in the setup of complex symmetrizable Kac-Moody algebras. More precisely, Littelmann defined root operators $E_p$ and $F_p$ on continuous paths $\pi\,:\,[0,1]\rightarrow \hR$ with $\pi(0)=0$. In fact, the operator $E_p$ (respectively $F_p$) is defined as in Remark \ref{defrootops} if the function $h_p\,:\,[0,1]\rightarrow {\mathbb R}$ given by $t\mapsto \langle\pi(t),\alpha_p^\vee\rangle$ is weakly decreasing (respectively weakly increasing) between $t_0$ and $t_1$. In general, the definition is more involved; for convenience, we stated it in Section \ref{sec:finite}. However, only the simpler version of the definition is needed for an LS chain, since it is known that the corresponding function $h_p$ satisfies the condition stated above. Littelmann considered the collection ${\mathcal P}_\lambda$ of all paths obtained by applying the operators $F_p$ to a fixed continuous path from $0$ to $\lambda$ which lies inside the dominant Weyl chamber. He showed that these paths form a crystal, that the associated colored directed graph does not depend on the initial path, and that one can express 
\begin{equation}\label{litchar}\chi(\lambda)=\sum_{\pi\in{\mathcal P}_\lambda} e^{\pi(1)}\,;\end{equation}
moreover, there is a corresponding Littlewood-Richardson rule \cite{Li1, Li2,Li3}. Stembridge \cite{St} reproved the special case of the above results corresponding to LS chains by showing that they form an admissible system. Kashiwara \cite{kascbm}, Lakshmibai \cite{lakbqd}, and Joseph \cite{josqgp} proved independently that Littelmann paths (obtained from a fixed path via root operators) have the structure of a perfect crystal.

Let us now return to LS chains, and fix $\lambda$ in $\Lambda^+$. Recall the set $I$ in (\ref{indset}), and the $\lambda$-chain $\{\beta_i\}_{i\in I}$ given by Proposition~\ref{constr-lambdachain}, which depends on a total order on the set of simple roots $\alpha_1<\cdots < \alpha_r$. We will now describe a bijection between the corresponding admissible subsets (cf.\ Definition~\ref{admsubset}) and the LS chains corresponding to the antidominant weight $-\lambda$. 

Given an index $i=(\alpha,k)$, we let $\beta_i:=\alpha$ and $t_i:=k/\langle\lambda,\alpha^\vee\rangle$. We have an order-preserving map from $I$ to $[0,1)$ given by $i\mapsto t_i$. Recall the notation $r_i:=s_{\beta_i}$ and $\widehat{r}_i:=s_{{\beta_i},l_i^\emptyset}$; we also let $\widehat{r}_i':=s_{{\beta_i},-l_i^\emptyset}$. Consider an admissible subset $J=\{j_1<j_2<\ldots<j_s\}$ and let
\[\{0=a_0<a_1<\ldots<a_l\}:=\{t_{j_1}\le t_{j_2}\le \ldots\le t_{j_s}\}\cup\{0\}\,.\]
Let $0=n_{0}\le n_1<\ldots<n_{l+1}=s$ be such that $t_{j_h}=a_k$ if and only if $n_{k}<h\le n_{k+1}$, for $k=0,\ldots,l$. Define Weyl group elements $u_h$ for $h=0,\ldots,s$ and $w_k$ for $k=0,\ldots,l$ by $u_0:=1$, $u_h:=r_{j_1}\ldots r_{j_h}$, and $w_k:=u_{n_{k+1}}$. Let also $\mu_k:=w_k(\lambda)$. For any $k=1,\ldots,l$, we have the following saturated chain in Bruhat order of minimum (left) coset representatives modulo $W_\lambda$:
\[w_{k-1}=u_{n_{k}}\lessdot u_{n_{k}+1}\lessdot\ldots\lessdot u_{n_{k+1}}=w_k\,;\]
indeed, none of the reflections $r_{j_1},\ldots,r_{j_s}$ lies in $W_\lambda$, since $\langle\lambda,\beta_i^\vee\rangle\ne 0$ for all $i\in I$. The above chain gives rise to a saturated increasing chain from $-\mu_{k-1}$ to $-\mu_k$ in the Bruhat order on $-W\lambda$. We will now show that this chain is, in fact, a chain in $a_k$-Bruhat order. Let $\widetilde{\beta}_h:=u_{h-1}(\beta_{j_h})$, so that $u_h:=s_{\widetilde{\beta}_h}u_{h-1}$, for $h=1,\ldots,s$. We need to check that $a_k\langle u_{h-1}(\lambda),\widetilde{\beta}_h^\vee\rangle\in{\mathbb Z}$, for $n_{k}<h\le n_{k+1}$. But
\[  \langle u_{h-1}(\lambda),\widetilde{\beta}_h^\vee\rangle=\langle\lambda,\beta_{j_h}^\vee\rangle\,,\]
while, by definition, $a_k=t_{j_h}$ is a fraction with denominator $\langle\lambda,\beta_{j_h}^\vee\rangle$. Hence 
\[-\mu_0<_{a_1}-\mu_1<_{a_2}\ldots<_{a_l}-\mu_l\]
is an LS chain in the $W$-orbit of $-\lambda$. We denote it by $\gamma(J)$, and the associated continuous piecewise-linear path by $\pi(J)$. 

Note that $\gamma(\emptyset)$ is the LS chain consisting only of $-\lambda$, while $\pi(\emptyset)\,:\,[0,1]\rightarrow\hR$ is the path $t\mapsto -t\,\lambda$. The path $\pi(\emptyset)$ intersects the affine hyperplane $H_{\beta_i,-l_i^\emptyset}$ at $t=t_i$ for $i\in I$; moreover, these and $t=1$ are the only intersections of $\pi(\emptyset)$ with affine hyperplanes $H_{\alpha,k}$, for $\alpha\in\Phi$ and $k\in {\mathbb Z}$. It is not hard to see that the path $\pi(J)$ can be described using folding operators as follows:
\begin{equation}\label{describefold}\pi(J):=\phi_{j_1}\ldots\phi_{j_s}(\pi(\emptyset))\,;\end{equation}
here, the folding operators $\phi_i$ are defined as follows on the relevant paths $\pi$:
\begin{equation}\label{simpledeffold}\phi_i(\pi)(t):=\casetwoex{\pi(t)}{0\le t\le t_i}{\widehat{r}_i'(\pi(t))}{t_i<t\le 1}\end{equation}

\begin{remark}\label{weightls} Based on (\ref{describefold}), we can easily show that $\pi(J)(1)=-\mu(J)$ (cf.\ Definition~{\rm \ref{defweighted}}). Indeed, we have $\widehat{r}_{j_1}\ldots \widehat{r}_{j_s}(\lambda)=-\widehat{r}_{j_1}'\ldots \widehat{r}_{j_s}'(-\lambda)$. \end{remark}

\begin{theorem}\label{rootopscorr}
We have 
\[E_p(\pi(J))=\pi(F_p(J))\]
for all admissible subsets $J$  (here $E_p$ is the root operator for paths, while $F_p$ in the one defined in Section {\rm \ref{sec:rootoperators}}).
\end{theorem}

\begin{proof}
Consider the point $P_\varepsilon:=\varepsilon\, \omega_1 +
\varepsilon^2 \omega_2 +\cdots +\varepsilon^r \omega_r$, 
where $\varepsilon$ is a small positive real number.
Let $\pi_{\varepsilon}:[0,1]\to \h_\R^*$ 
be the path $t\mapsto -t\,\lambda + P_\varepsilon$. Given $i=(\alpha,k)$ in $I$ and a sufficiently small $\varepsilon$, the path $\pi_\varepsilon$ crosses the affine hyperplane $H_{\alpha,-k}$ 
at $t_{\varepsilon,i}:=(k+\sum_{p=1}^r(\omega_p,\alpha^\vee)\,\varepsilon^p)/\langle\lambda,\alpha^\vee\rangle$. 
Fix a simple root $\alpha_p$ and an admissible subset $J=\{j_1<j_2<\ldots<j_s\}$. Let $\Gamma=\Gamma(J)$ be the corresponding admissible folding, and $L(\Gamma)=(l_i)_{i\in I}$ the corresponding level sequence. We can find $\varepsilon_*<1$ such that, for all $i<j$ in $J\cup I_{\alpha_p}(\Gamma)$ and $\varepsilon<\varepsilon_*$, the points $t_{\varepsilon,i},\,t_{\varepsilon,j}$ exist, and we have $t_{\varepsilon,i}<t_{\varepsilon,j}$. Let us now extend this path by the segments from 0 to $P_\varepsilon$ and from $-\lambda+P_\varepsilon$ to $-\lambda$. More precisely, we consider the path $\widehat{\pi}_{\varepsilon}:[0,1]\to \h_\R^*$ given by
\[\widehat{\pi}_\varepsilon(t):=\casethree{tP_\varepsilon}{0\le t<\varepsilon}{\pi_\varepsilon\left(\frac{t-\varepsilon}{1-2\varepsilon}\right)}{\varepsilon\le t\le 1-\varepsilon}{-\lambda+\frac{1-t}{\varepsilon}P_\varepsilon}{1-\varepsilon<t\le 1}\]
Since $\varepsilon_*<1$, the point $\widehat{\pi}_\varepsilon(t)$ does not lie on an affine hyperplane $H_{\alpha_p,k}$, for $k\in{\mathbb Z}$, whenever $t\in(0,\varepsilon]\cup[1-\varepsilon,1)$ and $\varepsilon<\varepsilon_*$. From now on, we assume that  $\varepsilon<\varepsilon_*$. Let $\widehat{t}_{\varepsilon,i}$ be the value of $t$ in $[\varepsilon,\,1-\varepsilon]$ for which $\widehat{\pi}_\varepsilon$ crosses the affine hyperplane $H_{\alpha,-k}$ for $i=(\alpha,k)$ in $I$ (assuming that such a value exists). Clearly, it is still true that, for all $i<j$ in $J\cup I_{\alpha_p}(\Gamma)$, the points $\widehat{t}_{\varepsilon,i},\,\widehat{t}_{\varepsilon,j}$ exist, and we have $\widehat{t}_{\varepsilon,i}<\widehat{t}_{\varepsilon,j}$. 

Now consider the path
\[\widehat{\pi}_\varepsilon(J):=\phi_{j_1}\ldots\phi_{j_s}(\widehat{\pi}_\varepsilon)\,,\]
where the folding operators $\phi_i$ are defined as in (\ref{simpledeffold}), except that $t_i$ is replaced by $\widehat{t}_{\varepsilon,i}$. It is easy to see that the only intersections of $\widehat{\pi}_\varepsilon(J)$ with affine hyperplanes $H_{\alpha_p,k}$, for $k\in{\mathbb Z}$, occur at $t=0$, $t=1$, and $t=\widehat{t}_{\varepsilon,i}$ for $i\in I_{\alpha_p}(\Gamma)$. Moreover, note that the signs of the roots in the pairs $(\delta_i,\zeta_i\delta_i)$ for $(\delta_i,\zeta_i)$ in $\Gamma$ with $\delta_i=\pm\alpha_p$, as well as the sign of $\langle\delta_\infty,\alpha_p^\vee\rangle$ indicate the sides of the mentioned hyperplane on which the path $\widehat{\pi}_\varepsilon$ lies before and/or after the intersection. Corollary \ref{cor1} shows that the minimum of the function $t\mapsto\langle\widehat{\pi}_\varepsilon(J)(t),\alpha_p^\vee\rangle$ is $-M(\Gamma,p)$; furthermore, this minimum is attained at $t=\widehat{t}_{\varepsilon,i}$ for $i\in I_{\alpha_p}(\Gamma)$ with $l_i=M(\Gamma,p)$, as well as at $t=1$ if $l_{\alpha_p}^\infty=M(\Gamma,p)$. Let $k:=k_F(\Gamma,p)$, $m:=m_F(\Gamma,p)$, and assume $m\ne\infty$. We define $a<b$ by $j_a<k<j_{a+1}$ and $j_b=m$ (possibly $a=0$, in which case the corresponding inequality is dropped). Let us define $E_p(\widehat{\pi}_\varepsilon(J))$ as in Remark \ref{defrootops}. By Corollary \ref{cor2}, the corresponding points $t_0$ and $t_1$ are as follows:
\[t_0=\widehat{t}_{\varepsilon,k}\,,\;\;\;\;t_1=\widehat{t}_{\varepsilon,m}\,.\]
Consider $i\in\{k<j_{a+1}<\ldots<j_{b-1}\}$, and let $i'$ be its successor in $J\cup\{k\}$. Given a subset $A=\{a_1<a_2<\ldots<a_l\}$ of $I$, we will use the notation $r_A$ for $r_{a_1}\ldots r_{a_l}$; we also set $r_\emptyset:=1$. Let $J':=\{j\in J\mid j< j_{a+1}\}$ and $J'':=\{j\in J\mid j_{a+1}\le j\le i\}$. By (\ref{defeppaths}), the direction of $E_p(\widehat{\pi}_\varepsilon(J))$ for $t\in[\widehat{t}_{\varepsilon,i},\,\widehat{t}_{\varepsilon,i'}]$ is $s_pr_{J'}r_{J''}(-\lambda)$. Since $F_p(J)=(J\cup\{k\})\setminus\{m\}$, the direction of $\widehat{\pi}_\varepsilon(F_p(J))$ for the same values of $t$ is $r_{J'}r_kr_{J''}(-\lambda)$. But $r_{J'}(\beta_k)=\alpha_p$, so $s_p=r_{J'}r_kr_{J'}^{-1}$, and, therefore, the two directions above coincide. The directions of the two paths also coincide for $t\ge \widehat{t}_{\varepsilon,m}$, since $\kappa(F_p(J))=\kappa(J)$ (cf.\ Proposition~\ref{closed}). The case $m=\infty$ is similar. All this shows that
\begin{equation}\label{foldepsilon}
E_p(\widehat{\pi}_\varepsilon(J))=\widehat{\pi}_\varepsilon(F_p(J))\,.
\end{equation}
Now let us take the limit as $\varepsilon\rightarrow 0$. Note that $\widehat{\pi}_\varepsilon(J)$ converges uniformly to $\pi(J)$, since $\widehat{t}_{\varepsilon,j_i}$ converges to $t_{j_i}$, for $i=1,\ldots,s$; so the minimum of $\pi(J)$ is also $-M(\Gamma,p)$, and the points $t_0,\,t_1$ in the construction of $E_p(\pi(J))$ are precisely $t_k,\,t_m$ (which are the limits of $\widehat{t}_{\varepsilon,k},\,\widehat{t}_{\varepsilon,m}$). This implies that $E_p$ commutes with limits, so the proposition follows by taking the limit in (\ref{foldepsilon}). 
\end{proof}

\begin{corollary}
The map $J\mapsto \gamma(J)$ is a bijection between the admissible subsets considered above and the LS chains corresponding to the antidominant weight $-\lambda$.
\end{corollary}

\begin{proof} Surjectivity follows directly from Theorem~\ref{rootopscorr}, based on the fact that all LS chains corresponding to $-\lambda$ can be obtained from the one consisting only of $-\lambda$ by applying the root operators $F_p$. Injectivity then follows from the character formula (\ref{litchar}), since $\pi(J)(1)=-\mu(J)$, as noted in Remark \ref{weightls}.
\end{proof}

\begin{remarks} (1) The proof of Theorem~\ref{rootopscorr} contains the justification of the fact that the minima of the paths associated to LS chains are integers. This justification is based only on the combinatorics in Section \ref{sec:folded}. Note that the same fact was proved by Littelmann in \cite{Li1} using different methods.

(2) The proof of Theorem~\ref{rootopscorr} shows that LS chains can be viewed as a limiting case of a special case of our construction. The special choices of $\lambda$-chains that lead to LS chains represent a very small fraction of all possible choices.
\end{remarks} 

Based on the independent results of Kashiwara \cite{kascbm}, Lakshmibai \cite{lakbqd}, and Joseph \cite{josqgp}, which were discussed above, we deduce the following corollary. 

\begin{corollary} Given a complex symmetrizable Kac-Moody algebra $\mathfrak g$, consider the colored directed graph defined by the action of root operators (cf.\ Section {\rm \ref{sec:rootoperators}}) on the admissible subsets corresponding to the special choice of a $\lambda$-chain above. This graph is isomorphic to the crystal graph of the irreducible representation with highest weight $\lambda$ of the associated quantum group $U_q({\mathfrak g})$.
\end{corollary}

We make the following conjecture, which is the analog of a result due to Littelmann, that was discussed above.

\begin{conjecture} The colored directed graph defined by the action of root operators on the admissible subsets corresponding to any $\lambda$-chain does not depend on the choice of this chain.
\end{conjecture}

The conjecture was proved for finite types (i.e., complex semisimple Lie algebras) in \cite{Le}. It implies that any choice of a $\lambda$-chain leads to a perfect crystal. 

\section{The Finite Case}
\label{sec:finite}

In this section, we discuss the way in which the model in this paper specializes to the one in \cite{LP} in the case of finite irreducible root systems.

Let $\Phi$ be the root system of a simple Lie algebra. Let $\Waff$ be  the {\it affine Weyl group\/} for $\Phi^\vee$, that is, the group generated by the affine reflections 
$s_{\alpha,k}$ (defined in (\ref{affref})). The corresponding affine hyperplanes $H_{\alpha,k}$ divide the real vector space $\hR$ into open
regions, called {\it alcoves.} The {\it fundamental alcove\/} $A_\circ$ is given by 
$$
A_\circ :=\{\lambda\in \hR \mid 0< \langle \lambda,\alpha^\vee \rangle <1 \textrm{ for all }
\alpha\in\Phi^+\}.
$$
We say that two alcoves are adjacent if they are
distinct and have a common wall.  
For a pair of adjacent alcoves, let us write 
$A\stackrel{\alpha}\longrightarrow B$ if the common wall of $A$ and $B$ 
is orthogonal to the root $\alpha\in\Phi$, and $\alpha$ points 
in the direction from $A$ to $B$.  

\begin{definition}
An {\it alcove path\/} is a sequence of alcoves
$(A_0,A_1,\dots,A_l)$ such that $A_{j-1}$ and $A_j$ are adjacent, for
$j=1,\dots,l$.
We say that an alcove path is {\it reduced\/} if it has minimal 
length among all alcove paths from $A_0$ to $A_l$.
\end{definition}

Let $A_\lambda=A_\circ + \lambda$ be the alcove obtained via the affine
translation of the fundamental alcove $A_\circ$ by a weight $\lambda$.  The 
reduced alcove paths from $A_\circ$ to $A_{\lambda}$ are in bijection with the
reduced decompositions of the element $v_\lambda$ in $\Waff$ defined by
$v_\lambda(A_\circ) = A_\lambda$; see~\cite{LP}.  Let us fix a dominant
weight $\lambda$.

\begin{proposition}  The sequence of roots 
$\{\beta_i\}_{i\in I}$  with $I=\{1,\ldots,l\}$ is a $\lambda$-chain 
(cf.\ Definition~{\rm \ref{lambdachain}}) if and only if   
 there exists a reduced alcove path
$A_0=A_\circ\stackrel{-\beta_1}\longrightarrow \cdots
\stackrel{-\beta_l}\longrightarrow A_{l}=A_{-\lambda}$.
\label{prop:alcoves=chains}
\end{proposition}

This proposition is an analog of the fact that the normal ordering of roots can be described
in terms of dihedral subgroups.  Each alcove $A$ is given by the inequalities
$$
A=\{v\in V\mid n_\alpha< \langle v,\alpha^\vee \rangle < n_\alpha+1,
\textrm{ for all roots }\alpha\in\Phi^+\},
$$
where $n_\alpha=n_\alpha(A)$ are some integers.  We need the following 
characterization due to Shi 
of the collection of integers $\{n_\alpha\}_{\alpha\in\Phi^+}$ associated with alcoves.

\begin{proposition} \cite{Shi}  An arbitrary collection of integers $\{m_\alpha\}_{\alpha\in\Phi^+}$
corresponds to some alcove $A$, i.e., $m_\alpha=n_\alpha(A)$ for all $\alpha\in\Phi^+$,
if and only if, for any triple of roots $\alpha,\beta,\gamma\in\Phi^+$ such that $\gamma^\vee=
\alpha^\vee+\beta^\vee$, we have $m_\gamma - m_\alpha-m_\beta\in\{0,1\}$.
\end{proposition}

\begin{proof}[Proof of Proposition~{\rm \ref{prop:alcoves=chains}}]
For a sequence of positive roots $(\beta_1,\dots,\beta_l)$,
define $m_\alpha^i := -\#\{j\leq i \mid \beta_j = \alpha\}$,
for $\alpha\in\Phi^+$ and $i=0,\dots,l$.  
The axioms of a $\lambda$-chain in the finite case (condition (1) in Definition~\ref{lambdachain} and condition ($2'$) in Proposition~\ref{equivalence}) can be rewritten
in terms of the integers $m_\alpha^i$ as follows:
(1) $0=m_\alpha^0\geq m_\alpha^1\geq \cdots \geq m_\alpha^l = 
- \langle \lambda,\alpha^\vee\rangle$, for $\alpha\in\Phi^+$; and
(2) for any triple $\alpha,\beta,\gamma\in\Phi^+$ such that $\alpha^\vee +\beta^\vee = \gamma^\vee$
and $i=0,\dots,l$, we have
$m_\gamma^i - m_\alpha^i-m_\beta^i \in\{0,1\}$ (interlacing condition).
Shi's result implies that these conditions are equivalent to the fact
that $A_0=A_\circ\stackrel{-\beta_1}\longrightarrow \cdots
\stackrel{-\beta_l}\longrightarrow A_{l}=A_{-\lambda}$ is a reduced alcove path,
where $A_i$ is the alcove associated with the collection of integers 
$\{m_\alpha^i\}_{\alpha\in\Phi^+}$, i.e., $n_\alpha(A_i)=m_\alpha^i$. 
\end{proof}

\begin{remarks} (1) In \cite{LP}, (reduced) $\lambda$-chains were defined as chains of roots determined by a reduced alcove path. As we have seen, the mentioned definition is equivalent to the one in this paper.

(2) Reduced alcove paths from $A_\circ$ to $w_\circ(A_\circ)=-A_\circ$ (where $w_\circ$ is the longest Weyl group element) correspond to {\em reflection orderings} \cite{Dyer}. If $\lambda$ is regular, then some reduced alcove paths from $A_\circ$ to $A_{-\lambda}$ start with an alcove path from $A_\circ$ to $w_\circ(A_\circ)$. Hence, we can say that $\lambda$-chains extend the notion of a reflection ordering. 
\end{remarks}

\begin{definition} A {\it gallery\/} is a sequence 
$\gamma=(F_0,A_0,F_1,A_1, F_2, \dots , F_l, A_l, F_{l+1})$ 
such that $A_0,\dots,A_l$ are alcoves;
$F_j$ is a codimension one common face of the alcoves $A_{j-1}$ and $A_j$,
for $j=1,\dots,l$; $F_0$ is a vertex of the first alcove $A_0$; and 
$F_{l+1}$ is a vertex of the last alcove $A_l$. 
Furthermore, we require that $F_0=\{0\}$, $A_0=A_\circ$, and 
$F_{l+1}=\{\mu\}$ for some weight $\mu\in\Lambda$,
which is called the {\it weight\/} of the gallery. 
The folding operator $\phi_j$ is the operator which acts on a gallery by  leaving its initial segment from $A_0$ to $A_{j-1}$ intact and by reflecting the remaining tail in the affine hyperplane containing the face $F_j$. In other words, we define
$$\phi_j(\gamma):=(F_0,A_0, F_1, A_1, \dots, A_{j-1}, F_j'=F_j,  A_{j}', F_{j+1}', A_{j+1}', \dots,  A_l', F_{l+1}'),$$
where $F_j\subset H_{\alpha,k}$, $A_i' := s_{\alpha,k}(A_i)$, and $F_i':=s_{\alpha,k}(F_i)$, for $i=j,\dots,l+1$.
\end{definition}

 The galleries defined above are special cases of the generalized galleries in~\cite{GaLi}.

Let us fix a reduced alcove path $A_0=A_\circ\stackrel{-\beta_1}\longrightarrow \cdots
\stackrel{-\beta_l}\longrightarrow A_{l}=A_{-\lambda}$, which determines the $\lambda$-chain
$\{\beta_i\}_{i\in I}$ with $I:=\{1,\ldots,l\}$. The alcove path also determines an obvious gallery 
\[\gamma(\emptyset)=(F_0,A_0,F_1,\dots,F_l, A_l, F_{l+1})\]
 of weight $-\lambda$. We use the same notation as in Sections \ref{sec:lambdachains}-\ref{sec:rootoperators}. For instance, $r_i:=s_{\beta_i}$ and $\widehat{r}_i:=s_{{\beta_i},l_i^\emptyset}$. We also let $\widehat{r}_i'$ be the affine reflection in the hyperplane containing $F_i$. 

\begin{definition} 
Given an admissible subset $J= \{ j_1<\cdots< j_s\}\subseteq I$ (cf.\ Definition~{\rm \ref{admsubset}}), we define the gallery $\gamma(J)$ as $\phi_{j_1}\cdots \phi_{j_s} (\gamma(\emptyset))$, and call it an {\it admissible folding\/} of $\gamma(\emptyset)$. 
\end{definition}

\begin{remark} The weight of the gallery $\gamma(J)$ is $-\mu(J)$ (cf.\ Definition~{\rm \ref{defweighted}}). Indeed, we have $\widehat{r}_{j_1}\ldots \widehat{r}_{j_s}(\lambda)=-\widehat{r}_{j_1}'\ldots \widehat{r}_{j_s}'(-\lambda)$. Hence, the model in this paper specializes to the one in \cite{LP}, whose construction was based on the geometry of the generalized flag variety.
\end{remark}

Since we assumed that $\Phi$ is irreducible, there is 
a unique {\it highest coroot\/} $\theta^\vee\in\Phi^\vee$, i.e., a unique coroot that has 
maximal height. We will also use the {\it Coxeter number}, that 
can be defined as $\hvee:=\langle\rho,\theta^\vee\rangle+1$ (in the finite case, the dominant weight $\rho$ considered at the end of Section \ref{sec:notation} is unique, and is given by $\frac{1}{2}\sum_{\alpha\in\Phi^+}\alpha$). Let $Z$ be the set of the elements of the lattice $\Lambda/\hvee$ 
that do not belong to any affine hyperplane $H_{\alpha,k}$. Each alcove $A$ contains precisely one element $\zeta_A$ of the set $Z$ (cf.\ \cite{Kost,LP}); this will be called the {\em central point} of $A$. In particular, $\zeta_{A_\circ}=\rho/\hvee$. 

\begin{proposition}\cite{LP}\label{adjalcoves} For a pair of adjacent alcoves $A\stackrel{\alpha}\longrightarrow B$, we have $\zeta_B -\zeta_A = \alpha/\hvee$.
\end{proposition}  

Let us now associate to the gallery $\gamma(\emptyset)$ a continuous piecewise-linear path. Consider the points $\eta_{0}:=0$, $\eta_{2i+1}:=\zeta_{A_i}$ for $i=0,\ldots,l$, $\eta_{2i}:=\frac{1}{2}(\eta_{2i-1}+\eta_{2i+1})$ for $i=1,\ldots,l$, and $\eta_{2l+2}:=-\lambda$. Note that $\eta_{2i}$ lies on $F_i$ for $i=0,\ldots,l+1$. Let $\pi(\emptyset)$ be the piecewise-linear path obtained by joining $\eta_{0},\,\eta_1,\,\ldots,\eta_{2l+2}$. Given an admissible subset $J$, let $\eta_{0}'=0,\,\eta_1'=\rho/\hvee,\,\eta_2',\ldots,\eta_{2l+2}'=-\mu(J)$ be the points on the faces of the gallery $\gamma(J)$ that are obtained (in the obvious way) from $\eta_{0},\,\eta_1,\,\eta_2,\,\ldots,\eta_{2l+2}$ in the process of constructing $\gamma(J)$ from $\gamma(\emptyset)$ via folding operators. Clearly, $\eta_{2i+1}'$ are the central points of the corresponding alcoves in $\gamma(J)$, for $i=0,\ldots,l$. By joining $\eta_{0}',\,\eta_1',\,\ldots,\eta_{2l+2}'$, we obtain a piecewise-linear path that we call $\pi(J)$. Note that $\pi(J)$ can be described using folding operators, as in (\ref{describefold}), once these operators are appropriately defined. 

\begin{remark} The maps $J\mapsto\gamma(J)$ and $J\mapsto \pi(J)$ are one-to-one. Indeed, given the gallery $\gamma(J)=(F_0,A_0,F_1,\dots, A_l,F_{l+1})$, we have $J=\{j\mid A_{j-1}=A_j\}$. Also, the gallery $\gamma(J)$ can be easily recovered from $\pi(J)$.
\end{remark}

The next result follows from Proposition~\ref{adjalcoves} and the definition of folding operators on chains of roots and galleries. 

\begin{proposition}\label{chain_path} Let $\Gamma(J)=(\{(\gamma_i,\varepsilon_i)\}_{i\in I},\gamma_\infty)$. Then, for all $i\in I$, we have
\[{\eta_{2i-1}'-\eta_{2i}'}=\frac{\gamma_i}{2\hvee}\,,\;\;\;\;\;\;\;{\eta_{2i}'-\eta_{2i+1}'}=\frac{\varepsilon_i\gamma_i}{2\hvee}\,,\;\;\;\;\;\;\;{\eta_{2l+1}'-\eta_{2l+2}'}=\frac{\gamma_{\infty}}{\hvee}\,.\]
\end{proposition}

It turns out that, in general, the collection of paths $\pi(J)$, for $J$ ranging over admissible subsets, does {\em not} coincide with the collection of Littelmann paths obtained from $\pi(\emptyset)$ by applying the root operators $E_p$. Indeed, it is not true in general that $E_p(\pi(J))=\pi(F_p(J))$, as was the case with the paths corresponding to LS chains (cf.\ Theorem~\ref{rootopscorr}). The reason is that, given $\pi=\pi(J)$, the function $h_p\,:\,[0,1]\rightarrow {\mathbb R}$ given by $t\mapsto \langle\pi(t),\alpha_p^\vee\rangle$ is usually not weakly decreasing between the corresponding points $t_0$ and $t_1$ (see Section \ref{sec:ls} for the definition of these points). This happens, for instance, when applying $E_2$ to the path $\pi(\emptyset)$ in Example \ref{exg2} below. Such situations can also arise if we define $\pi(\emptyset)$ by joining the centers of the faces $F_i$, or the centers of both the alcoves $A_i$ and the faces $F_i$ (in the order they appear in the gallery $\gamma(\emptyset)$). In all these situations, we need the general definition of $E_p$ for Littelmann paths, which we now recall from \cite{Li2}.

As in Section \ref{sec:ls}, the definition is easier to state if we replace the continuous piecewise-linear paths $\pi$ (satisfying $\pi(0)=0$) with their left-hand derivatives $\gamma$ (which are piecewise-constant left-continuous maps defined on $(0,1]$). Recall that we denoted the minimum of the function $h_p$ by $m_p$. If $m_p>-1$, then $E_p(\gamma)$ is undefined, as before, so assume that $m_p\le -1$. Choose $t_0=x_0<x_1<\ldots<x_l=t_1$ such either
\begin{enumerate}
\item $h_p(x_{i-1})=h_p(x_i)$ and $h_p(t)\ge h_p(x_{i-1})$ for $t\in [x_{i-1},\,x_i]$;
\item or $h_p$ is strictly decreasing on $[x_{i-1},\,x_i]$ and $h_p(t)\ge h_p(x_{i-1})$ for $t\le x_{i-1}$.
\end{enumerate}
We now define
\begin{equation}\label{gendefeppaths}E_p(\gamma)(t):=\casetwo{s_p(\gamma(t))}{\mbox{$x_{i-1}<t\le x_i$ and $h_p$ behaves on $[x_{i-1},\,x_i]$ as in (2)}}{\gamma(t)}\end{equation}

\begin{remark} Our model can be based on {\em any} $\lambda$-chain (that is, not necessarily on the ones given by Proposition~\ref{constr-lambdachain}), and we still have a definition of root operators for admissible subsets/admissible foldings that corresponds to the simpler version of their definition on paths, given in (\ref{defeppaths}).
\end{remark}

\begin{figure}[ht]
\input{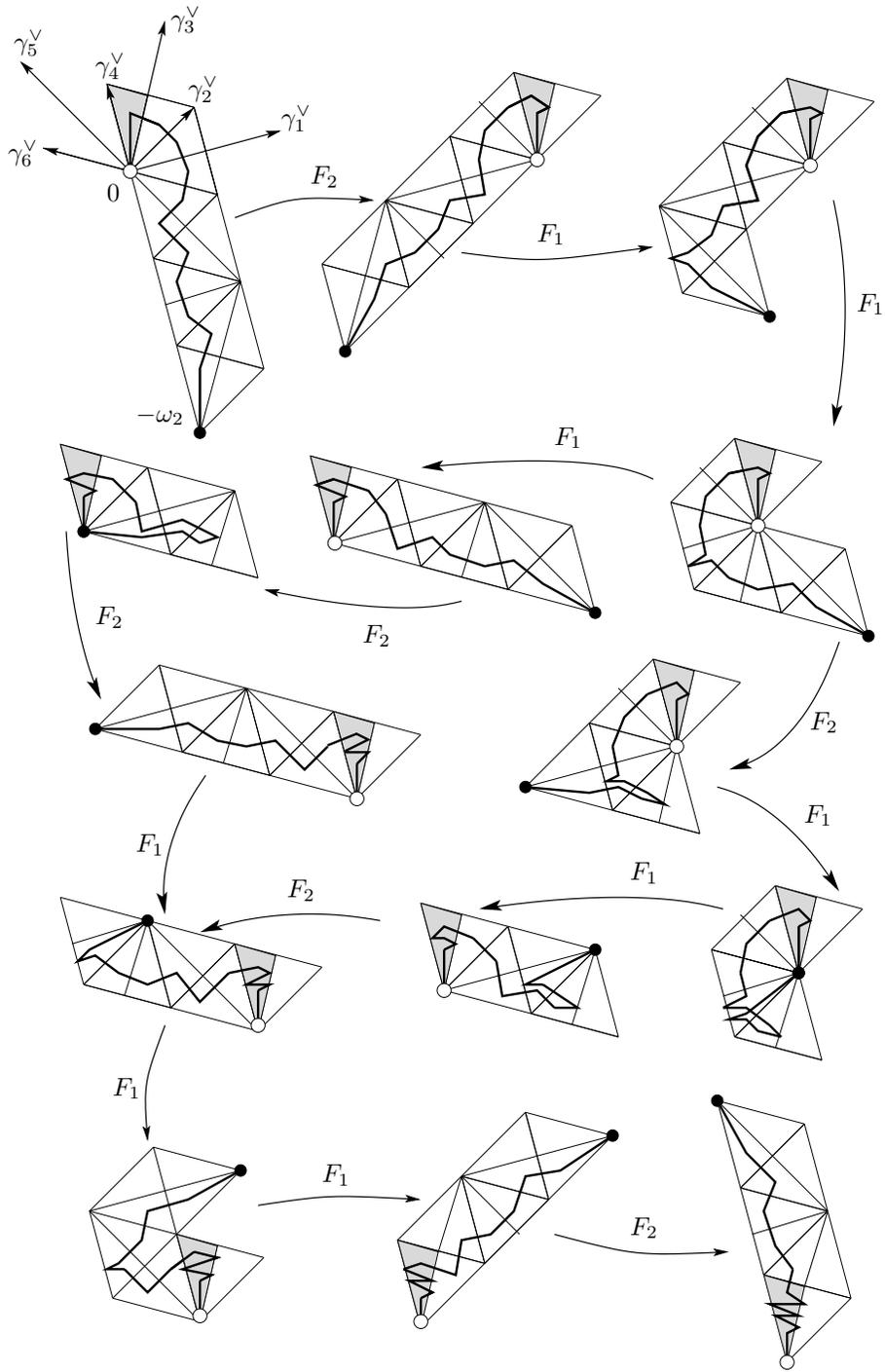}
\caption{The crystal for the fundamental weight $\omega_2$ for type $G_2$.} 
\label{fig:G2}
\end{figure}

\begin{example}\label{exg2}  Suppose that the root system $\Phi$ is of type $G_2$.
The positive roots are $\gamma_1=\alpha_1,\ 
\gamma_2=3\alpha_1+\alpha_2,\ \gamma_3=2\alpha_1+\alpha_2,\ 
\gamma_4=3\alpha_1+2\alpha_2,\ \gamma_5=\alpha_1+\alpha_2,\
\gamma_6=\alpha_2$.
The corresponding coroots are
$\gamma_1^\vee=\alpha_1^\vee, \
\gamma_2^\vee= \alpha_1^\vee+\alpha_2^\vee,\
\gamma_3^\vee = 2\alpha_1^\vee+3\alpha_2^\vee,\
\gamma_4^\vee = \alpha_1^\vee+2\alpha_2^\vee,\
\gamma_5^\vee = \alpha_1^\vee+3\alpha_2^\vee,\ 
\gamma_6^\vee=\alpha_2^\vee$.

Suppose that $\lambda=\omega_2$.
Proposition~\ref{constr-lambdachain} gives the following $\omega_2$-chain:
$$
(\beta_1,\dots,\beta_{10})=
({\gamma_6},{\gamma_5},{\gamma_4},{\gamma_3},{\gamma_2},{\gamma_5},{\gamma_3},
{\gamma_4}, {\gamma_5},{\gamma_3})\,.
$$
Thus, we have $\widehat{r}_1=s_{\gamma_6,0}$, $\widehat{r}_2=s_{\gamma_5,0}$,
$\widehat{r}_3=s_{\gamma_4,0}$, $\widehat{r}_4=s_{\gamma_3,0}$,
$\widehat{r}_5=s_{\gamma_2,0}$, $\widehat{r}_6=s_{\gamma_5,1}$,
$\widehat{r}_7=s_{\gamma_3,1}$, $\widehat{r}_8=s_{\gamma_4,1}$,
$\widehat{r}_9=s_{\gamma_5,2}$, $\widehat{r}_{10}=s_{\gamma_3,2}$. There are
six saturated chains in the Bruhat order (starting at the identity) on the
corresponding Weyl group that can be retrieved as subchains of the
$\omega_2$-chain. We indicate each such chain and the corresponding admissible
subsets in $\{1,\dots,10\}$.
\begin{enumerate}
\item 1: $\{\}$;
\item $1<s_{\gamma_6}$: $\{1\}$;
\item $1<s_{\gamma_6}<s_{\gamma_6}s_{\gamma_5}$: $\{1,2\}$, $\{1,6\}$, $\{1,9\}$;
\item $1<s_{\gamma_6}<s_{\gamma_6}s_{\gamma_5}<s_{\gamma_6}s_{\gamma_5}s_{\gamma_4}$: $\{1,2,3\}$, $\{1,2,8\}$, $\{1,6,8\}$;
\item $1<s_{\gamma_6}<s_{\gamma_6}s_{\gamma_5}<s_{\gamma_6}s_{\gamma_5}s_{\gamma_4}<s_{\gamma_6}s_{\gamma_5}s_{\gamma_4}s_{\gamma_3}$: $\{1,2,3,4\}$, $\{1,2,3,7\}$, $\{1,2,3,10\}$, $\{1,2,8,10\}$, $\{1,6,8,10\}$;
\item $1<s_{\gamma_6}<s_{\gamma_6}s_{\gamma_5}<s_{\gamma_6}s_{\gamma_5}s_{\gamma_4}<s_{\gamma_6}s_{\gamma_5}s_{\gamma_4}s_{\gamma_3}<s_{\gamma_6}s_{\gamma_5}s_{\gamma_4}s_{\gamma_3}s_{\gamma_2}$: $\{1,2,3,4,5\}$.
\end{enumerate}
The weight of each admissible subset is now easy to compute (by applying the
corresponding affine reflections above to $\omega_2$, cf.\ Definition
\ref{defweighted}). This leads to the expression for the character $\chi(\omega_2)$
as the following sum over admissible subsets: 
$$
\chi(\omega_2) = e^{\omega_2}  + e^{\widehat{r}_1 (\omega_2)}+
e^{\widehat{r}_1\, \widehat{r}_2(\omega_2)}+
e^{\widehat{r}_1\, \widehat{r}_6(\omega_2)}+
e^{\widehat{r}_1\, \widehat{r}_9(\omega_2)}+
\cdots +
e^{\widehat{r}_1 \,\widehat{r}_6 \,\widehat{r}_8\, \widehat{r}_{10}(\omega_2)} +
e^{\widehat{r}_1 \, \widehat{r}_2 \, \widehat{r}_3 \, \widehat{r}_4 \, \widehat{r}_5(\omega_2)}.
$$

Figure~\ref{fig:G2} displays the galleries $\gamma(J)$ corresponding to the
admissible subsets $J$ indicated above, the associated paths $\pi(J)$, as well
as the action of the root operators $F_p$ on $J$. 
For each path, we shade the fundamental alcove, mark the origin by a white dot 
``$\circ$'', and mark the endpoint of a black dot ``$\bullet$''.
Since some linear steps in $\pi(J)$ might coincide, we
display slight deformations of these paths, so that no information is lost in
their graphical representations.  As discussed above, the weights of the
irreducible representation $V_{\omega_2}$ are obtained by changing the signs of
the endpoints of the paths $\pi(J)$ (marked by black dots).
The roots in the corresponding admissible foldings $\Gamma(J)$ can also be read off;
see Proposition~\ref{chain_path}. 
At each step, a path $\pi(J)$ either crosses a wall of the affine Coxeter arrangement
or bounces off a wall.  The associated admissible subset $J$
is the set of indices of bouncing steps in the path.

\label{example:G-2}
\end{example}


\phantom{a}\bigskip

\begin{thebibliography}{Dyer}

\bibitem[Ch]{Chev}
C.~Chevalley.
\newblock Sur les d\mbox{\'{e}}compositions cellulaires des espaces
  \mbox{$G/B$}.
\newblock In {\em Algebraic \mbox{G}roups and \mbox{G}eneralizations:
  \mbox{C}lassical \mbox{M}ethods}, volume 56 Part 1 of {\em Proceedings and
  Symposia in Pure Mathematics}, pages 1--23. Amer.\ Math.\ Soc., 1994.

\bibitem[De]{Deo1}
V.~V. Deodhar.
\newblock Some characterizations of {B}ruhat ordering on a {C}oxeter group and
  determination of the relative {M}\"obius function.
\newblock {\em Invent.\ Math.}, 39:187--198, 1977.

\bibitem[Dyer]{Dyer}
M.~J. Dyer.
\newblock Hecke algebras and shellings of {B}ruhat intervals.
\newblock {\em Compositio Math.}, 89(1):91--115, 1993.

\bibitem[Fu]{fulyt}
W.~Fulton.
\newblock {\em Young Tableaux}, volume~35 of {\em London Math. Soc. Student
  Texts}.
\newblock Cambridge Univ. Press, Cambridge and New York, 1997.

\bibitem[GL]{GaLi}
S.~Gaussent and P.~Littelmann.
\newblock {LS}-galleries, the path model and {MV}-cycles.
\newblock {\em Duke Math. J.},  127:35--88, 2005.

\bibitem[Hu]{Hum}
J.~E. Humphreys.
\newblock {\em Reflection Groups and {C}oxeter Groups}, volume~29 of {\em
  Cambridge Studies in Advanced Mathematics}.
\newblock Cambridge University Press, Cambridge, 1990.

\bibitem[Jos]{josqgp}
A.~Joseph.
\newblock {\em Quantum Groups and Their Primitive Ideals}, 
\newblock Springer-Verlag, New York, 1994.

\bibitem[Kac]{kacidl}
V.~G. Kac.
\newblock {\em Infinite Dimensional {L}ie Algebras}.
\newblock Cambridge University Press, Cambridge, 1990.

\bibitem[Ka1]{kascqa}
M.~Kashiwara.
\newblock Crystalizing the {$q$}-analogue of universal enveloping algebras.
\newblock {\em Commun. Math. Phys.}, 133:249--260, 1990.

\bibitem[Ka2]{kascbq}
M.~Kashiwara.
\newblock On crystal bases of the {$q$}-analogue of universal enveloping
  algebras.
\newblock {\em Duke Math. J.}, 63:465--516, 1991.

\bibitem[Ka3]{kascbm}
M.~Kashiwara.
\newblock Crystal bases of modified quantized enveloping algebra.
\newblock {\em Duke Math.\ J.}, 73:383--413, 1994.

\bibitem[Kos]{Kost}
B.~Kostant.
\newblock Powers of the {E}uler product and commutative subalgebras of a
  complex simple {L}ie algebra.
\newblock{\em Invent. Math.}, 158:181--226, 2004.

\bibitem[Ku]{kumkmg}
S.~Kumar.
\newblock {\em Kac-{M}oody Groups, Their Flag Varieties and Representation
  Theory}, volume 204 of {\em Progress in Mathematics}.
\newblock Birkh\"auser Boston Inc., Boston, MA, 2002.

\bibitem[La]{lakbqd}
V.~Lakshmibai.
\newblock Bases for quantum {D}emazure modules.
\newblock In {\em Representations of Groups (Banff, AB, 1994)}, volume~16 of
  {\em CMS Conf. Proc.}, pages 199--216. Amer.\ Math.\ Soc., Providence, RI,
  1995.

\bibitem[LS]{LS1}
V.~Lakshmibai and C.~S. Seshadri.
\newblock Standard monomial theory.
\newblock In {\em Proceedings of the Hyderabad Conference on Algebraic Groups
  (Hyderabad, 1989)}, pages 279--322, Madras, 1991. Manoj Prakashan.

\bibitem[LSc]{lasksb}
A.~Lascoux and M.-P. Sch\mbox{\"{u}}tzenberger.
\newblock Keys and standard bases.
\newblock In D.~Stanton, editor, {\em Invariant Theory and Tableaux}, volume~19
  of {\em The IMA Vol. in Math. and Its Appl.}, pages 125--144,
  Berlin-Heidelberg-New York, 1990. Springer-Verlag.

\bibitem[Le]{Le}
C.~Lenart.
\newblock On the combinatorics of crystal graphs, I. Lusztig's involution.
\newblock {\tt arXiv:math.RT/0509200}, to appear in {\em Adv. Math.}

\bibitem[LP]{LP}
C.~Lenart and A.~Postnikov.
\newblock Affine {W}eyl groups in {$K$}-theory and representation theory.
\newblock {\tt arXiv:math.RT/0309207}.

\bibitem[Li1]{Li1}
P.~Littelmann.
\newblock {A Littlewood-Richardson rule for symmetrizable Kac-Moody algebras}.
\newblock {\em Invent.\ Math.}, 116:329--346, 1994.

\bibitem[Li2]{Li2}
P.~Littelmann.
\newblock Paths and root operators in representation theory.
\newblock {\em Ann.\ of Math.\ (2)}, 142:499--525, 1995.

\bibitem[Li3]{Li3}
P.~Littelmann.
\newblock Characters of representations and paths in $\hR$.
\newblock {\em Proc.\ Sympos.\ Pure Math.}, 61:29--49, 1997.

\bibitem[Lu]{luscba}
G.~Lusztig.
\newblock Canonical bases arising from quantized enveloping algebras. {II}.
\newblock {\em Progr. Theoret. Phys. Suppl.}, (102):175--201, 1991.

\bibitem[PR]{PR}
H.~Pittie and A.~Ram.
\newblock A \mbox{P}ieri-\mbox{C}hevalley formula in the \mbox{$K$}-theory of a
  \mbox{$G/B$}-bundle.
\newblock {\em Electron.\ Res.\ Announc.\ Amer.\ Math.\ Soc.}, 5:102--107, 1999.

\bibitem[Shi]{Shi} 
J.-Y.~Shi.
\newblock Sign type corresponding to an affine Weyl group.
\newblock {\em J.\ London Math.\ Soc.\ (2)} 35:56--74, 1987.


\bibitem[Ste]{St}
J.~R. Stembridge.
\newblock Combinatorial models for {W}eyl characters.
\newblock {\em Adv.\ Math.}, 168:96--131, 2002.

\end{thebibliography}

\end{document}